\documentclass[12pt]{article}
\usepackage{amsmath,amsfonts}
\usepackage{amssymb}
\usepackage{amsthm}
\usepackage{times,txfonts}

  \usepackage{kbordermatrix}

\usepackage{pict2e}

\newtheorem{thm}{Theorem}[section]
\newtheorem{prop}[thm]{Proposition}
\newtheorem{cor}[thm]{Corollary}
\newtheorem{lemma}[thm]{Lemma}
\newtheorem{dfn}[thm]{Definition}

\newtheorem{remark}[thm]{\it Remark}

\numberwithin{equation}{section}

\def\pf{\noindent{\it Proof.} \ }

\def\qed{\hfill $\square$}

\def\id{{\rm id}}

\title{Hermite-Pad\'e approximation, isomonodromic deformation and hypergeometric integral}

\author{Toshiyuki Mano and Teruhisa Tsuda}
\date{February 24, 2015
\\
Revised: March 29, 2016 (final version)}

\begin{document}
\maketitle

\renewcommand{\thefootnote}{\fnsymbol{footnote}}
\footnotetext{ 
 (secondary).  
} 

\footnotetext{ \ \\
Toshiyuki Mano  \quad  {\tt tmano@math.u-ryukyu.ac.jp}
\\
Department of Mathematical Science,  University of the Ryukyus,
Okinawa 903-0213, Japan
\\ \\
Teruhisa Tsuda  \quad  {\tt tudateru@econ.hit-u.ac.jp}
\\
Department of Economics, Hitotsubashi University,  
Tokyo 186-8601, Japan}

\begin{abstract}
We develop
an underlying relationship between the theory of rational approximations 
and that of isomonodromic deformations.
We show that
a certain duality in Hermite's two approximation problems 
for functions
leads to the Schlesinger transformations,
i.e.  transformations of a  linear differential equation
shifting its characteristic exponents by integers while keeping its monodromy 
invariant.
Since approximants 
and remainders are described by block-Toeplitz determinants,
one can clearly understand the determinantal structure in isomonodromic deformations.
We demonstrate our method in a certain family of Hamiltonian systems of isomonodromy type 
including the sixth Painlev\'e equation and Garnier systems;
particularly,  
we present their
solutions written in terms of iterated hypergeometric integrals.
An algorithm for constructing the Schlesinger transformations is also discussed through vector continued fractions.
\end{abstract}

\tableofcontents

\section*{Introduction}

Let $L$ be an integer greater than one.
For a given $L$-tuple 
of analytic functions
(or formal power series) 
$f_0(w), f_1(w), \ldots, f_{L-1}(w)$
and of nonnegative integers
${\boldsymbol m}=(m_0,m_1, \ldots,m_{L-1})$,
Hermite 
considered
the following two rational approximation problems.
The first is to find 
$L$ polynomials
\[
{\mathfrak q}^{{\boldsymbol m}}_i(w) \quad (0 \leq i \leq L-1)
\]
of degree 
$m_i-1$
such that 
\[
\sum_{i=0}^{L-1}
f_i {\mathfrak q}^{{\boldsymbol m}}_i
=O(w^{ |{\boldsymbol m}| -1}),
\]
where
$|{\boldsymbol m}|= \sum_{i=0}^{L-1}m_i$;
i.e., 
the left-hand side has a zero of order at least
$|{\boldsymbol m}| -1$
at $w=0$.
The second is to find $L$ polynomials
\[
\mathfrak{p}^{{\boldsymbol m}}_i(w) \quad (0 \leq i \leq L-1)
\]
of degree 
$|{\boldsymbol m}|-m_i$ 
such that 
\[
f_i \mathfrak{p}^{{\boldsymbol m}}_j-f_j \mathfrak{p}^{{\boldsymbol m}}_i=O(w^{ |{\boldsymbol m}| +1}).
\]
Each system of polynomials 
$\{ {\mathfrak q}^{{\boldsymbol m}}_i\}$ and $\{ {\mathfrak p}^{{\boldsymbol m}}_i\}$
generically
turns out to be  
unique up to simultaneous multiplication by constants,
as an elementary consequence of linear algebra.
The above approximation problems come from Hermite's study on 
arithmetic properties of the exponential function
and 
are called collectively  
the {\it Hermite-Pad\'e approximations};
often,
the former 
is referred to 
as the `type I' and the latter 
as the `type II' or as the {\it simultaneous Pad\'e approximation}.
Note that if $L=2$ then both of them reduce to the
(usual) Pad\'e approximations.
Although these two types of approximations were seemingly unrelated, 
Mahler discovered that they were fundamentally 
connected to each other.
Put 
${\boldsymbol e}_i=(0, \ldots,0,
\stackrel{ \stackrel{i}{\smile} }{1},0,\ldots,0) \in {\mathbb Z}^L$
for $0 \leq i \leq L-1$.

\begin{thm}[Mahler's duality]
\label{thm:original}
It holds that
\[
\begin{bmatrix}
{\mathfrak q}^{{\boldsymbol m}+{\boldsymbol e}_0 }_0 &{\mathfrak q}^{{\boldsymbol m}+{\boldsymbol e}_0 }_1 & \cdots&{\mathfrak q}^{{\boldsymbol m}+{\boldsymbol e}_0 }_{L-1} 
\\
{\mathfrak q}^{{\boldsymbol m}+{\boldsymbol e}_1 }_0 &{\mathfrak q}^{{\boldsymbol m}+{\boldsymbol e}_1 }_1 & \cdots&{\mathfrak q}^{{\boldsymbol m}+{\boldsymbol e}_1 }_{L-1} 
\\
\vdots & \vdots& \ddots & \vdots
\\
{\mathfrak q}^{{\boldsymbol m}+{\boldsymbol e}_{L-1} }_0 &{\mathfrak q}^{{\boldsymbol m}+{\boldsymbol e}_{L-1} }_1 & \cdots&{\mathfrak q}^{{\boldsymbol m}+{\boldsymbol e}_{L-1} }_{L-1} 
\end{bmatrix}
\begin{bmatrix}
{\mathfrak p}^{{\boldsymbol m}-{\boldsymbol e}_0 }_0 &{\mathfrak p}^{{\boldsymbol m}-{\boldsymbol e}_1 }_0 & \cdots&{\mathfrak p}^{{\boldsymbol m}-{\boldsymbol e}_{L-1} }_{0} 
\\
{\mathfrak p}^{{\boldsymbol m}-{\boldsymbol e}_0 }_1 &{\mathfrak p}^{{\boldsymbol m}-{\boldsymbol e}_1 }_1 & \cdots&{\mathfrak p}^{{\boldsymbol m}-{\boldsymbol e}_{L-1} }_{1} 
\\
\vdots & \vdots& \ddots & \vdots
\\
{\mathfrak p}^{{\boldsymbol m}-{\boldsymbol e}_{0} }_{L-1} &{\mathfrak p}^{{\boldsymbol m}-{\boldsymbol e}_{1} }_{L-1} & \cdots&{\mathfrak p}^{{\boldsymbol m}-{\boldsymbol e}_{L-1} }_{L-1} 
\end{bmatrix}
=w^{|{\boldsymbol m}|} \cdot D
\]
with $D$ being a diagonal constant matrix.
Moreover, if every diagonal part
${\mathfrak q}^{{\boldsymbol m}+{\boldsymbol e}_i }_i$
and 
${\mathfrak p}^{{\boldsymbol m}-{\boldsymbol e}_i }_i$
 is chosen to be a monic polynomial then $D$ becomes the identity matrix.
\end{thm}

The aim of this paper is to develop an underlying relationship of 
Hermite's two approximations
with
the theory of linear differential equations in the complex domain,
especially with that of isomonodromic deformations.
Interestingly enough, 
Mahler's duality plays a crucial role in constructing
a certain class of 
{\it Schlesinger transformations},
i.e.
transformations of a linear differential equation
shifting its characteristic exponents by integers while keeping its monodromy invariant.

In Sect.~\ref{sect:hpa},
we begin by introducing 
the two types of rational approximations for an $L$-tuple of functions,
which are slightly modified
(in order to fit the construction of Schlesinger transformations)
from the original Hermite-Pad\'e and simultaneous Pad\'e approximations.
We then prove a variation of Mahler's duality between them (see Theorem~\ref{thm:dual}).
Applying the approximations for the solution of an $L \times L$ Fuchsian system of linear differential equations
yields its Schlesinger transformation
(see Theorem~\ref{thm:schl});
in fact,
Mahler's duality guarantees
 the absence of apparent singularities in the new Fuchsian system
after the Schlesinger transformation.
In Sect.~\ref{sect:det},
we deduce determinantal representations for the approximants and remainders
 from the approximation conditions
(see Propositions~\ref{prop:detq}  and \ref{prop:detp}
 and also Remark~\ref{remark:bT}).
 It should be noted that any key ingredient here 
is
 written in terms of {\it block-Toeplitz determinants}.
In Sect.~\ref{sect:cf},
we present an algorithm for constructing  the Schlesinger transformation via vector continued fraction expansions,
which is  a variation of the {\it Jacobi--Perron algorithm}
(i.e. a higher dimensional analogue of the Euclidean algorithm).

The last two sections are devoted to the study of isomonodromic deformations.
In Sect.~\ref{sect:aid},
we first review  the {\it Schlesinger system}
of nonlinear differential equations,
which
governs  isomonodromic deformations of a Fuchsian system.
Since a Schlesinger transformation preserves
 the monodromy of the
Fuchsian system, 
it gives rise to a discrete symmetry of the corresponding Schlesinger system.
We clarify,
based on the above relationship with rational approximations,
 the determinantal structure in the general solutions of the Schlesinger systems.
Next we concern
a particular family
of the Schlesinger systems,
which possesses a unified description as a
polynomial Hamiltonian system denoted by ${\cal H}_{L,N}$ ($L \geq 2, N \geq 1$);
it includes various noteworthy examples 
of isomonodromic deformations 
such as 
the sixth Painlev\'e equation (${\cal H}_{2,1}$)
and the Garnier system in $N$ variables (${\cal H}_{2,N}$).
In  Sect.~\ref{sect:ihi},
we demonstrate 
Schlesinger transformations on the previously known hypergeometric solution of ${\cal H}_{L,N}$ (see \cite{tsu12b});
as a result,
we obtain
solutions of ${\cal H}_{L,N}$ 
written in terms of iterated hypergeometric integrals
through Fubini's theorem and the Vandermonde determinant
(see Theorem~\ref{thm:hgi} and its sequel).

\section{From Hermite's two approximation problems to Schlesinger transformations}
\label{sect:hpa}

In this section 
we show how rational approximations are
useful for constructing
Schlesinger transformations
of linear differential equations.
Fix an integer $L \geq 2$.
We shall first introduce 
two different types of rational approximation
problems for 
an $L$-tuple   
\[f_0(w), f_1(w), \ldots, f_{L-1}(w)  \in {\mathbb C}[\![w]\!]\]
of formal power series,
where we assume $f_0(0) \neq 0$ without loss of generality.

\subsection{Hermite-Pad\'e approximation (of type I)}
\label{subsect:hpa}

Let $n$ be a positive integer. 
Consider for each $i$ $(0 \leq i \leq L-1)$ an $L$-tuple of polynomials
$Q_j^{(i)}=Q_j^{(i)}(w)$ $(0 \leq j \leq L-1)$
of degree at most $n-1+\delta_{i,j}$,
where
\[
\delta_{i,j} =
\left\{ 
\begin{array}{ll}
1& (i = j)
\\
0 & (i \neq j)
\end{array}
\right.
\]
is the Kronecker delta.
Suppose
the approximation condition
\begin{equation} \label{eq:hpa}
 Q_0^{(i)} f_0+\cdots +Q_i^{(i)} f_i +w Q_{i+1}^{(i)} f_{i+1}  + \cdots + w Q_{L-1}^{(i)}f_{L-1} = O(w^{nL})
 \end{equation}
 is  fulfilled for each $i$.
This condition
amounts to
 a system of  $n L$  homogeneous linear equations
for the $n L+1$
unknown coefficients of the polynomials $Q_j^{(i)}$ $(0 \leq j \leq L-1)$.
Under a generic condition for the power series
$f_0, \ldots,f_{L-1}$,
these polynomials  are uniquely determined up to simultaneous multiplication by constants;
see Sect.~\ref{subsect:dethpa}.
We will be concerned with the row vector
\[\widetilde{\boldsymbol Q}^{(i)}
=\left( \widetilde{Q}^{(i)}_j \right)_{0 \leq j \leq L-1}
 =\left( Q_0^{(i)}, \ldots ,Q_i^{(i)}, w Q_{i+1}^{(i)},  \ldots ,w Q_{L-1}^{(i)} \right).
\]

\begin{remark} \rm
By construction,
the polynomial
$Q^{(0)}_0(w)$ 
has no constant term.
Moreover,
the degree of the {\it diagonal} part  $Q^{(i)}_i(w)$ $(0 \leq i \leq L-1)$
 turns out to be $n$ exactly;
see (\ref{eq:deg}) in Sect.~\ref{subsect:dual}.
\end{remark}

\subsection{Simultaneous Pad\'e approximation}
\label{subsect:spa}

We treat another type of approximation problem
for the same power series $f_0, \ldots,f_{L-1}$.
Consider for each $j$ $(0 \leq j \leq L-1)$ an $L$-tuple of polynomials 
$ P_i^{(j)}=P_i^{(j)}(w)$ $(0 \leq i \leq L-1)$
of degree at most $n (L-1)-1+\delta_{i,j}$
which satisfies
the following approximation conditions:
\begin{equation}  \label{eq:spa}
\begin{array}{llll}
\bullet \ \text{if $j=0$} &&
f_0 P_i^{(0)}-f_i  P_0^{(0)}=O(w^{nL}) & \text{for $1 \leq i\leq L-1$};
\\ 
\bullet \ \text{if $1 \leq  j \leq L-1$} &&
f_0 P_i^{(j)}-f_i  P_0^{(j)}=O(w^{nL}) & \text{for $1 \leq i \leq j-1$},
\\
&&
f_0 P_i^{(j)}-f_i w P_0^{(j)}=O(w^{nL}) & \text{for $j \leq i \leq L-1$}.
\end{array}
\end{equation}
These conditions
are interpreted as 
a system of  $n L(L-1)$  homogeneous linear equations
for the $n L(L-1)+1$
unknown coefficients of the polynomials $P_i^{(j)}$
$(0 \leq i \leq L-1)$.
Hence 
the column vector
\begin{align*}
 \widetilde{\boldsymbol P}^{(j)}
 ={}^{\rm T}\left( \widetilde{P}_i^{(j)} \right)_{0 \leq i \leq L-1}
={}^{\rm T}\left( w P_0^{(j)}, \ldots , w P_{j-1}^{(j)}, P_{j}^{(j)},  \ldots ,P_{L-1}^{(j)} \right)
\end{align*}
is generically unique up to multiplication by constants;
see Sect.~\ref{subsect:detspa}.

\begin{remark}\rm
Let $a \neq b$.
It is immediate from $1/f_0 \in {\mathbb C}[\![w] \!]$
to verify that
  \[
 \begin{array}{lll}
 \bullet \ \text{\rm if $a, b <j$} &&
 f_a \widetilde{P}_b^{(j)} - f_b \widetilde{P}_a^{(j)} =O(w^{nL+1});
 \\
  \bullet \ \text{\rm otherwise} &&
  f_a \widetilde{P}_b^{(j)} - f_b\widetilde{P}_a^{(j)} =O(w^{nL}).
 \end{array}
 \]
\end{remark}

\subsection{Mahler's duality in the two approximation problems}
\label{subsect:dual}

There is an interesting connection between the two approximation problems 
(\ref{eq:hpa}) and (\ref{eq:spa})
although they are seemingly unrelated.
The following theorem is thought of as a variation of 
{\it Mahler's duality}; 
see Theorem~\ref{thm:original} or \cite{mah68}.
We will give a proof of it
because our setup is 
slightly 
different from the original;
cf. \cite[Theorem~8.1.2]{bg-m96}
and
\cite{coa66}.

\begin{thm}
\label{thm:dual}
 It holds that
\[
\begin{bmatrix} 
\widetilde{\boldsymbol Q}^{(0)}
\\
\vdots
\\
\widetilde{\boldsymbol Q}^{(L-1)}
\end{bmatrix}
\left[
\widetilde{\boldsymbol P}^{(0)}, \ldots,\widetilde{\boldsymbol P}^{(L-1)}
\right]=w^{nL} \cdot D
\]
with $D$ being a diagonal constant matrix.
\end{thm}

\pf
Let us first estimate the degree of the $(i,j)$-entry
\[
M_{ij}=\sum_{k=0}^{L-1} \widetilde{Q}_k^{(i)} \widetilde{P}_k^{(j)}
\]
of the left-hand side.
The degree of each polynomial $\widetilde{Q}_k^{(i)}$ or $\widetilde{P}_k^{(j)}$ 
reads 
\begin{align*}
&
\underbrace{Q_0^{(i)}, \ldots ,Q_{i-1}^{(i)}}_{\text{$\deg \leq n-1$}},  
\underbrace{Q_i^{(i)}, w Q_{i+1}^{(i)},  \ldots ,w Q_{L-1}^{(i)}}_{\text{$\deg \leq n$ }},
\\ 
&
\underbrace{w P_0^{(j)}, \ldots , w P_{j-1}^{(j)}, P_{j}^{(j)}}_{\text{$\deg \leq n(L-1)$ }},
\underbrace{P_{j+1}^{(j)}, \ldots ,P_{L-1}^{(j)} }_{\text{$\deg \leq n(L-1)-1$ }}.
\end{align*}
Hence we have
\begin{equation}
 \label{eq:degM}
\deg_w M_{ij} \leq
\left\{ 
\begin{array}{ll}
n L & (i \leq j)
\\
n L-1 & (i > j).
\end{array}
\right.
\end{equation}

Next we shall estimate the multiplicity of $M_{ij}$
at $w=0$
by means of the
 approximation conditions.
Consider 
\begin{align*}
f_0  M_{ij} &= f_0 \sum_{k=0}^{L-1} \widetilde{Q}_k^{(i)} \widetilde{P}_k^{(j)}
\\
&= f_0\widetilde{Q}_0^{(i)} \widetilde{P}_0^{(j)}
+ \sum_{k \neq 0} \widetilde{Q}_k^{(i)} f_0 \widetilde{P}_k^{(j)}.
\end{align*}

\paragraph{}
\underline{(i) Case $i<j$ (strictly upper triangular part)}
\begin{align*}
\sum_{k \neq 0} \widetilde{Q}_k^{(i)} f_0 \widetilde{P}_k^{(j)}
&= \sum_{0 <  k <j  }
\widetilde{Q}_k^{(i)} f_0\widetilde{P}_k^{(j)}
+
\sum_{j \leq k}
\widetilde{Q}_k^{(i)} f_0 \widetilde{P}_k^{(j)}
\\
&=
\sum_{0 < k <j  }
\widetilde{Q}_k^{(i)} \left(f_k \widetilde{P}_0^{(j)} +O(w^{nL+1}) \right)
+
\sum_{j \leq k}
\widetilde{Q}_k^{(i)} \left( f_k \widetilde{P}_0^{(j)}+O(w^{nL}) \right),
\qquad 
\text{using (\ref{eq:spa})}
\\
&= \widetilde{P}_0^{(j)} \sum_{k \neq 0} \widetilde{Q}_k^{(i)} f_k 
+O(w^{nL+1}),
\qquad \text{since $\widetilde{Q}_k^{(i)}$ is divisible by $w$ if $i<(j \leq ) k$.}
\end{align*}
Therefore,
\begin{align*}
f_0 M_{ij}
&=  \widetilde{P}_0^{(j)} \sum_k \widetilde{Q}_k^{(i)} f_k +O(w^{nL+1})
\\
&=  \widetilde{P}_0^{(j)} \cdot O(w^{nL}) +O(w^{nL+1}),
\qquad 
\text{using (\ref{eq:hpa})}
\\
&= O(w^{nL+1}), \qquad 
\text{
since $\widetilde{P}_0^{(j)}$ is divisible by $w$ if $0(\leq i) <j$.}
\end{align*}

\paragraph{}
\underline{(ii) Case $i \geq j$ (lower triangular part)}
\begin{align*}
\sum_{k \neq 0} \widetilde{Q}_k^{(i)} f_0 \widetilde{P}_k^{(j)}
&= \sum_{0<k<j}  \widetilde{Q}_k^{(i)} f_0 \widetilde{P}_k^{(j)}
+  \sum_{j \leq k}  \widetilde{Q}_k^{(i)} f_0 \widetilde{P}_k^{(j)}
\\
&=\sum_{0 < k <j  }
\widetilde{Q}_k^{(i)} \left(f_k \widetilde{P}_0^{(j)} +O(w^{nL+1}) \right)
+
\sum_{j \leq k}
\widetilde{Q}_k^{(i)} \left( f_k \widetilde{P}_0^{(j)}+O(w^{nL}) \right),
\qquad 
\text{using (\ref{eq:spa})}
\\ 
&= \widetilde{P}_0^{(j)}\sum_{k \neq 0} \widetilde{Q}_k^{(i)} 
f_k +O(w^{nL}).
\end{align*}
Therefore,
\begin{align*}
f_0 M_{ij}
&=
 \widetilde{P}_0^{(j)} \sum_{k} \widetilde{Q}_k^{(i)} f_k +O(w^{nL})
 \\
 &=\widetilde{P}_0^{(j)} \cdot O(w^{nL})+O(w^{nL}),
 \qquad 
\text{using (\ref{eq:hpa})}
\\
&=O(w^{nL}).
\end{align*}

Noticing $1/f_0 \in {\mathbb C}[\![w]\!]$ we verify 
\[
M_{ij} = \left\{ 
\begin{array}{ll}
O(w^{nL+1}) & (i <j)
\\
O(w^{nL})& (i \geq j).
\end{array}
\right.
\]
Combining  this with (\ref{eq:degM}),
we can conclude that
$(M_{ij}) = w^{nL} \cdot D$.
\qed
\\

Consequently, the diagonal entry $M_{ii}$
coincides with the term of highest degree in $Q_i^{(i)}(w)P_i^{(i)}(w)$
and thus
\begin{equation}
\label{eq:deg}
\deg_w Q_i^{(i)}=n, \quad \deg_w P_i^{(i)}=n(L-1).
\end{equation}
We henceforth normalize $\widetilde{\boldsymbol Q}^{(i)}$ and $\widetilde{\boldsymbol P}^{(j)}$
so that their diagonal parts 
 $Q^{(i)}_i$ and $P^{(j)}_j$ become
  monic polynomials,
 i.e.
 \[
\left. z^{n}Q^{(i)}_i(z^{-1}) \right|_{z=0}=
\left.z^{n(L-1)}P^{(i)}_i(z^{-1}) \right|_{z=0}=1
\]
and thereby $D= I$ (the identity matrix).

\begin{cor} \label{cor:R}

The polynomial matrix 
\begin{equation} \label{eq:smatrix}
R(z)
=
z^n
\begin{bmatrix} 
\widetilde{\boldsymbol Q}^{(0)}(z^{-1})
\\
\vdots
\\
\widetilde{\boldsymbol Q}^{(L-1)}(z^{-1})
\end{bmatrix}
\in {\mathbb C}^{L \times L}[z]
\end{equation}
satisfies
\[
\begin{array}{ll}
{\rm(i)} & 
R(z)^{-1}
=z^{n(L-1)}
\left[
\widetilde{\boldsymbol P}^{(0)}(z^{-1}), \ldots,
 \widetilde{\boldsymbol P}^{(L-1)}(z^{-1})\right];
 \\
 {\rm(ii)} &  \det R(z)=1,
 \quad
\text{ i.e. }R(z) \in {\rm SL}(L, {\mathbb C}[z]).
\end{array}
\]
\end{cor}

\pf
Theorem~\ref{thm:dual}
shows (i) immediately.
Then, it holds that $R \in {\rm GL}(L, {\mathbb C}[z])$
and thus $\det R \in {\mathbb C} \setminus \{0\}$.
By definition,
$R$  
takes the form
\begin{equation} \label{eq:R&R-1}
R=\begin{bmatrix}
1 &  \cdots &*\\  
&\ddots&  \vdots \\
&&1
\end{bmatrix} 
+O(z);
\end{equation}
namely, 
its constant term is an upper triangular matrix whose 
diagonal entries are all one.
Therefore, we have $\det R =1$.
\qed

\subsection{Schlesinger transformations}
\label{subsect:st}

Consider an $L\times L$ Fuchsian system
\begin{equation} \label{eq:fs}
\frac{{\rm d} Y}{{\rm d} z}
=A Y
=\sum_{i=0}^{N+1} \frac{A_i}{z-u_i} Y
\quad
(A_i \text{ : constant matrix})
\end{equation}
of linear ordinary
differential equations
with 
$N+3$ regular singularities
\[
S=\{u_0=1, u_1,\ldots, u_N,u_{N+1}=0, u_{N+2}=\infty\}
\subset 
{\mathbb P}^1={\mathbb C} \cup \{\infty\}
\]
on the Riemann sphere.
Let
$A_{N+1}$ and $A_{N+2}=-\sum_{i=0}^{N+1} A_i$
be upper and lower triangular matrices, respectively.
Assume  
for simplicity
there is no integer difference among the 
{\it characteristic exponents} 
$\{{\varepsilon_{0,j}}\}_{0 \leq j \leq L-1}$ at $z=0$
(resp. $\{{\varepsilon_{\infty,j}}\}_{0 \leq j \leq L-1}$ at $z=\infty$),
i.e. the eigenvalues of the residue matrix $A_{N+1}$ (resp. $A_{N+2}$).
Then we have a solution 
$Y=Y(z)$ of (\ref{eq:fs}) normalized as 
\begin{align}
Y&=
\left(
\begin{bmatrix}
1 & \cdots &*\\  
&\ddots& \vdots \\
&&1
\end{bmatrix} 
+O(z)
\right)
\cdot
{\rm diag \ }  (z^{\varepsilon_{0,j}})_{0 \leq j \leq L-1}
\nonumber
\\ 
&=\Phi (w) \cdot {\rm diag \ }  (w^{\varepsilon_{\infty,j}})_{0 \leq j \leq L-1} \cdot C
\label{eq:fss}
\end{align}
with $w=1/z$ and
$C$ being an invertible constant matrix
(the {\it connection matrix} between $z=0$ and $z=\infty$).
Here $\Phi(w)$ is a matrix function holomorphic
at $w=0$ ($z=\infty$)
and $\Phi(0)$ is  invertible and lower triangular,
i.e.
\[
\Phi(w)
=
\begin{bmatrix}
* &  &\\  
 \vdots&\ddots& \\ 
*& \cdots&*
\end{bmatrix}
+O(w).
\]

\begin{remark}\rm
Many literatures adopt a different 
normalization such that the residue matrix at $z=\infty$ becomes diagonal.
Our present normalization
treats the two points $z=0$ and $z=\infty$ equally and it
emerges naturally from the similarity reduction of the UC hierarchy,
which is a context of infinite-dimensional integrable systems;
see \cite{tsu04, tsu12a, tsu14a}.
Furthermore, as clarified by Haraoka \cite{har12},
this normalization is effective to find a 
`good' coordinate of the space of Fuchsian systems having a given 
Riemann scheme.
\end{remark}

An analytic continuation along a loop on ${\mathbb P}^1 \setminus S$
based at some point $z_0$
induces a linear transformation of  $Y$
according to its multi-valuedness at the branch points $S$.
We thus obtain an
$L$-dimensional representation of the fundamental group
$\pi_1({\mathbb P}^1 \setminus S;z_0)$,
which is called the {\it monodromy} of $Y$. 
A left multiplication
$Y \mapsto \hat{Y}=RY$
of a rational  function matrix $R=R(z)$
is said to be a {\it Schlesinger transformation} if
the new equation
\[
\frac{{\rm d} \hat{Y}}{{\rm d} z}
=\hat{A} \hat{Y}
\]
satisfied by $\hat Y$ becomes the same form as the original (\ref{eq:fs}).
Because $R(z)$ is rational, 
$\hat Y$ and $Y$ have the same monodromy though they have different characteristic exponents by integers.
It is known that if we specify an admissible discrete change of the 
 characteristic exponents then
 the corresponding rational function matrix $R$
 of the Schlesinger transformation 
  is algebraically  computable from $Y$;
see \cite{jm81, sch12}.
In fact, the construction problem of Schlesinger transformations
is naturally related to rational approximation problems; see also Remark~\ref{remark:history}.

In this paper
we focus on
a class of Schlesinger transformations,
which is of particular interest from the viewpoint of 
Hermite's two approximation problems
and also of vector continued fractions (see Sect.~\ref{sect:cf}).
Note that, for a general Schlesinger transformation
other than the present direction,
though it can also be controlled by some rational approximation problems
but
it becomes much more complicated due to the absence 
of a duality like Mahler's; 
e.g.  $R^{-1}$ seems not to have a concise 
determinantal representation.

Let us define the $L$-tuple
${\boldsymbol f}={}^{\rm T}(f_0, \ldots,f_{L-1})$ of power series 
in $w$ 
as the first column of $\Phi(w)$,
where $\Phi(w)$ is the power series part of the solution $Y$ 
of (\ref{eq:fs}) near $z=\infty$.
Notice that $f_0(0)\neq 0$  certainly holds.
Therefore, all the general arguments in Sects.~\ref{subsect:hpa}--\ref{subsect:dual} are still valid for this specific case, and we are led to 
the Schlesinger transformation through the two approximation problems
for ${\boldsymbol f}$.
Now we state the result.

\begin{thm} 
 \label{thm:schl}
The polynomial matrix $R(z)$ given by  {\rm(\ref{eq:smatrix})}
realizes the Schlesinger transformation
$Y \mapsto \hat{Y}=R Y$ 
shifting 
the characteristic exponents at 
$z=\infty$
by
 $(n(L-1),-n,\ldots,-n)$.
\end{thm}

\pf 
It follows from the Hermite-Pad\'e approximation condition (\ref{eq:hpa})
that 
$R  {\boldsymbol f}=O(w^{n (L-1)})$.
By definition, $R$ takes the form 
\begin{equation} \label{eq:R}
R
= 
w^{-n}
\left(
\begin{bmatrix}
0 &  & &\\ 
*&*&  & \\ 
\vdots&\ddots&\ddots    &\\ 
*&\cdots&*&*
\end{bmatrix} 
+O(w)
\right). 
\end{equation}
Hence, if  we write as
\[
R \Phi(w) = 
\hat{\Phi}(w) 
\cdot {\rm diag \ } 
(w^{n(L-1)}, w^{-n}, \ldots,w^{-n}),
\]
then $\hat{\Phi}$ becomes the same form as $\Phi$.
Also, we verify from (\ref{eq:R&R-1}) that
$R$ does not change the form of the power series expansion of $Y$
near $z=0$.

On the other hand, 
the coefficient
\[
A(z)=
\sum_{i=0}^{N+1} \frac{A_i}{z-u_i} 
\]
of the Fuchsian system  (\ref{eq:fs}) is transformed as
\begin{equation}
\label{eq:AtoA}
A \mapsto \hat{A}=R A R^{-1} + \frac{{\rm d} R}{{\rm d}z}R^{-1} .
\end{equation}
If we remember both $R$ and $R^{-1}$ being polynomials
(see Corollary~\ref{cor:R}),
then $\hat{A}$ turns out to be
a rational function matrix having only simple poles at $S$ 
as well as the original $A$.
In this sense
Mahler's duality guarantees the absence of apparent singularities in the new equation ${\rm d} \hat{Y}/{{\rm d} z}
=  \hat{A} \hat{Y}$
satisfied by $ \hat{Y}=RY$.
\qed

\begin{remark}\rm
In the rank two case ($L=2$)
a similar construction of 
Schlesinger transformations
as Theorem~\ref{thm:schl}
has been established in 
\cite{man12}
based on 
 (usual)
Pad\'e approximations.
\end{remark}

\begin{remark} \rm
\label{remark:history}
A series of pioneering works
was done by  D. Chudnovsky and  G. Chudnovsky 
on
the close connection between 
rational approximation problems 
and Riemann's monodromy problem,
involving (semi-classical) orthogonal polynomials;
see \cite{cc82, cc94} and references therein. 
The `Pad\'e method'  recently
proposed by Yamada \cite{yam09}
is a recipe for Lax  formalism
of Painlev\'e equations and, at the same time,
for their special solutions,
which is based on
Pad\'e approximations (or interpolations) of
elementary functions;
interestingly enough,
it is applicable also for various discrete analogues of 
Painlev\'e equations beyond the originals;
see 
\cite{ika13, nag15, nty12}.

The essential idea of the above works could be exemplified by the  
following:
let us consider a function
$\varphi(z)= z^a(z-1)^b(z-u)^c$
with $a+b+c=0$.
The remainder 
$\rho:=P \varphi -Q =O(z^{-n-1})$
of its 
Pad\'e approximation
then satisfies
a second-order linear differential equation
denoted by $E$,
which may have 
an apparent singularity besides the
four regular singularities $S=\{0,1,u,\infty\} \subset {\mathbb P}^1$.
However, the two functions $\varphi$ and $\rho$
share the same multi-valuedness
since they are rationally related;
thus, the monodromy of $E$ 
is obviously constant with respect to $u$.
This fact
leads to special solutions of 
the sixth Painlev\'e equation $P_{\rm VI}$,
i.e. the {\it isomonodromic deformation}
(cf. Sect.~\ref{sect:aid})
 of a
second-order linear differential equation
with four regular singularities.

It is interesting to note that
such an idea had been recognized
implicitly by Laguerre
(before the discovery of Painlev\'e equations);
see \cite{lag80} and also \cite{mag95}.
\end{remark}

\begin{remark}\rm
The approximation conditions (\ref{eq:hpa}) and (\ref{eq:spa})
can be interpreted as certain
multi-orthogonality relations among the $L$-tuples of 
polynomials 
$\widetilde{\boldsymbol Q}^{(i)}$
and
$\widetilde{\boldsymbol P}^{(j)}$,
respectively;
i.e., these polynomials
can constitute multi-orthogonal polynomial systems.
In this paper,
although we do not enter into details on such aspects, 
we present below 
the {\it determinantal representations} for them,
which will crucially work in the last two sections.
\end{remark}

\section{Determinantal representations for approximation polynomials and remainders}

\label{sect:det}

In this section
we derive determinantal representations for 
the approximation polynomials $Q^{(i)}_j(w)$ and $P^{(i)}_j(w)$.
We
write the power series as
\begin{equation} \label{eq:fi}
f_i(w)=\sum_{j=0}^\infty a_j^iw^j \in {\mathbb C}[\![w]\!]
\end{equation}
henceforth;
note that the superscript $i$ of $a_j^i$
is just an index, not an exponent.
Introduce the $k \times l$ 
{\it rectangular  Toeplitz  matrix}
\begin{align}
A^i_j(k,l)
&=
\left[a^i_{j+m-n}\right]_{
\begin{subarray}{l}1 \leq m \leq k 
\\1 \leq n \leq l
\end{subarray}}
\nonumber
\\
&= 
\kbordermatrix{%
&1&2& &l  \\
1&a^i_j & a^i_{j-1} & \cdots &a^i_{j-l+1} \\
2&a^i_{j+1} & a^i_j & \cdots &a^i_{j-l+2} \\
  & \vdots& \vdots &              &\vdots   \\
k&a^i_{j+k-1} & a^i_{j+k-2} & \cdots &a^i_{j+k-l} 
}
\label{eq:rTm}
\end{align}
for the sequence $\{  a^i_j \}_{j=0}^\infty$,
where $a^i_j=0$ if $j<0$.
It holds that
\begin{equation}
\label{eq:tenchi}
\left(A^k_n(m,n)\right)_{i,j}
= \left( A^k_{m}(n,m) \right)_{n-j+1,m-i+1}
\end{equation}
by definition.

\subsection{Hermite-Pad\'e polynomials}
\label{subsect:dethpa}
We can calculate separately for each $i$ $(0 \leq i \leq L-1)$.
Therefore, for brevity,
we shall express the coefficients of
the approximation polynomial
$Q_j^{(i)}(w)$
as
\[
Q_j^{(i)}(w)=\sum_{k=0}^{ n-1+\delta_{i,j}  } b_{j,k}w^k
\]
with omitting the index $i$. 
The condition (\ref{eq:hpa})
implies vanishing of
the coefficients of 
$1, w, w^2, \ldots, w^{nL-1}$ in the left-hand side.
Consequently, we have a system 
\begin{equation}
\label{eq:lhpa}
{\cal A}^{(i)}\begin{bmatrix}
{\boldsymbol b}_0 \\
{\boldsymbol b}_1 \\
 \vdots \\{\boldsymbol b}_{L-1}
\end{bmatrix}
= 
\left.
 \begin{bmatrix}
0 \\ 0 \\\vdots \\ 0
\end{bmatrix}
\right\} \text{\footnotesize $nL$}
\end{equation}
of homogeneous linear equations
for the $nL+1$ unknowns
\[{\boldsymbol b}_j
={}^{\rm T}( b_{j,0},  \ldots , b_{j,n-1+\delta_{i,j}})
\quad
(0 \leq j \leq L-1)
\]
where 
\[
{\cal A}^{(i)}
=
\left[
\underbrace{A^0_0(nL,n)  \ \cdots \   A^{i-1}_0(nL,n) \
 \  A^i_0(nL,n+1)}_{\text{$i+1$ blocks}}
   \ 
\underbrace{A^{i+1}_{-1}(nL,n)  \
  \cdots    \
 A^{L-1}_{-1}(nL,n)}_{\text{$L-i-1$ blocks}}
\right].
\]
The solution of (\ref{eq:lhpa}) is unique up to multiplication by constants
if and only if  the rank of the
$nL \times (nL+1)$ matrix
${\cal A}^{(i)}$ equals $nL$
(which we will always assume).

Interestingly enough, we have the following determinantal representation of $Q^{(i)}_j(w)$.

\begin{prop}  \label{prop:detq}
It holds that
\begin{equation} \label{eq:Q}
Q^{(i)}_j(w) 
=\frac{1}{{\rm NQ}^{(i)}}
\det \left[
\begin{array}{ccc}
{\boldsymbol 0} &
\overbrace{1, w, \ldots, w^{n-1+\delta_{i,j}}}^{\text{\rm $j$th block}}
 &{\boldsymbol 0} 
 \\  \hline
&{\cal A}^{(i)}&
\end{array}
\right],
\end{equation}
where ${\rm NQ}^{(i)}$ are some normalizing constants. 
\end{prop}

\pf
Consider
\[
\rho_i(w)= Q_0^{(i)} f_0+\cdots +Q_i^{(i)} f_i +w Q_{i+1}^{(i)} f_{i+1}  + \cdots + w Q_{L-1}^{(i)}f_{L-1},
\]
which is  the remainder of the approximation condition (\ref{eq:hpa}).
Substituting (\ref{eq:Q}) shows that
\begin{align*}
\rho_i(w)
&= \frac{1}{{\rm NQ}^{(i)}}
\det
\left[
\begin{array}{cccc}
\overbrace{f_0, f_0 w,\ldots,f_0 w^{n-1}}^{\text{$0$th block}} 
& \cdots 
& \overbrace{f_{i-1}, f_{i-1} w,\ldots,f_{i-1} w^{n-1}}^{\text{$(i-1)$th block}}
& \overbrace{f_i, f_i w,\ldots,f_i w^n}^{\text{$i$th block}}
\\
A^0_0(nL,n) 
& \cdots 
&  A^{i-1}_0(nL,n) 
&  A^{i}_0(nL,n+1) 
\end{array}
\right.
\\
& \qquad \qquad \qquad \qquad \qquad
\left.
\begin{array}{ccc}
\overbrace{f_{i+1}w, f_{i+1} w^2,\ldots,f_{i+1} w^n}^{\text{$(i+1)$th block}}
& \cdots 
&\overbrace{f_{L-1}w, f_{L-1} w^2,\ldots,f_{L-1} w^n}^{\text{$(L-1)$th block}}
\\
A^{i+1}_{-1}(nL,n) 
& \cdots   
& A^{L-1}_{-1}(nL,n)
\end{array}
\right].
\end{align*}
Therefore, if we put $\rho_i(w)=\sum_{k=0}^\infty \rho_{k}^i w^k$, 
then the coefficients read 
\[
\rho_{k}^i
=
\frac{1}{{\rm NQ}^{(i)}}
\det
\begin{bmatrix}
A^0_k(1,n) 
&  \cdots 
&  A^{i-1}_k(1,n)
&  A^{i}_k(1,n+1)
& A^{i+1}_{k-1}(1,n)
& \cdots
& A^{L-1}_{k-1}(1,n)
\\
A^0_0(nL,n) 
& \cdots 
&  A^{i-1}_0(nL,n) 
&  A^{i}_0(nL,n+1) 
&A^{i+1}_{-1}(nL,n) 
& \cdots   
& A^{L-1}_{-1}(nL,n)
\end{bmatrix}.
\]
It is immediate
from a property of determinants
to verify 
$\rho_{k}^i=0$
for any $k$ less than $nL$;
thus, we have $\rho_i(w)=O(w^{nL})$ indeed.
\qed
\\

We will normalize the polynomials so that its diagonal part
$Q^{(i)}_i(w)$ becomes monic 
as well as in Sect.~\ref{subsect:dual}.
Accordingly, the normalizing constant ${\rm NQ}^{(i)}$ should be
\begin{align*}
{\rm NQ}^{(i)}
&=\det \left[
\begin{array}{ccc}
{\boldsymbol 0} &
\overbrace{0, \ldots, 0, 1}^{\text{\rm $i$th block}}
 &{\boldsymbol 0} 
 \\  \hline
&{\cal A}^{(i)}&
\end{array}
\right]
\\
&=(-1)^{n(i+1)} \det \begin{bmatrix}
A^0_0(nL,n) & \cdots 
&  A^{i}_{0}(nL,n) &A^{i+1}_{-1}(nL,n)& \cdots   
& A^{L-1}_{-1}(nL,n)
\end{bmatrix}.
\end{align*}
Thus,
the leading coefficient $\rho^i_{nL}$
of the remainder is given by
\begin{align*}
\rho^i_{nL}
&=
\frac{(-1)^{nL}}{{\rm NQ}^{(i)}}
\det
\left[
\begin{array}{cccc}
A^0_0(nL+1,n) 
& \cdots 
&  A^{i-1}_0(nL+1,n) 
&  A^{i}_0(nL+1,n+1) 
\end{array}
\right.
\\
& \qquad \qquad \qquad
\left.
\begin{array}{cccc}
&A^{i+1}_{-1}(nL+1,n) 
& \cdots   
& A^{L-1}_{-1}(nL+1,n)
\end{array}
\right];
\end{align*}
the constant term of the polynomial $Q_i^{(i)}(w)$
$(i \neq 0)$
is given by
\begin{align*}
Q_i^{(i)}(0)
&=\frac{1}{{\rm NQ}^{(i)}}
\det \left[
\begin{array}{ccc}
{\boldsymbol 0} &
\overbrace{1,0, \ldots, 0}^{\text{\rm $i$th block}}
 &{\boldsymbol 0} 
 \\  \hline
&{\cal A}^{(i)}&
\end{array}
\right]
\\
&=\frac{(-1)^{n i}}{{\rm NQ}^{(i)}} 
\det \begin{bmatrix}
A^0_0(nL,n) & \cdots 
&  A^{i-1}_{0}(nL,n) &A^{i}_{-1}(nL,n)& \cdots   
& A^{L-1}_{-1}(nL,n)
\end{bmatrix}.
\end{align*}

\begin{remark}\rm
\label{remark:bT}
We here restrict ourselves to the case where $f_0(w)=1$. 
Let us introduce the {\it block-Toeplitz determinant}
\[
\Delta^{(i)}
=
\det
 \left[
\underbrace{
A^1_{n}( m,  n) \
\cdots  \
A^{i-1}_{n}(m,  n)}_{\text{$i-1$ blocks}}
 \
\underbrace{
A^{i}_{n-1}(m,  n)  \
\cdots   \  
A^{L-1}_{n-1}(m,  n)
}_{\text{$L-i$ blocks}}
\right]
\]
of size $m=n(L-1)$
for each $i$ $(1 \leq i \leq L)$;
e.g.
${\rm NQ}^{(i)}
=(-1)^{n(i+1)} \Delta^{(i+1)}$.
We see in particular that 
\begin{align*}
\rho^0_{nL}
&=(-1)^m \frac{\Delta^{(L)}}{\Delta^{(1)}},
\\
Q_i^{(i)}(0)
&=(-1)^n
\frac{\Delta^{(i)}}{\Delta^{(i+1)}}
\quad (1 \leq i \leq L-1).
\end{align*}
These simple formulae will be used later in Sect.~\ref{sect:ihi}.
\end{remark}

\subsection{Simultaneous Pad\'e polynomials}
\label{subsect:detspa}

Suppose $f_0(w) = 1$ for simplicity.
Or, equivalently, we may understand that
we have renamed $f_i/f_0$ $(i \neq 0)$ as $f_i$.
For a given power series $F(w)=\sum_{k=0}^\infty  F_k w ^k$,
we employ the notation
\[
\left[  F(w)  \right]^{b}_{a}
=\sum_{k=a}^b  F_k w ^k
\]
denoting its section between $w^a$ and $w^b$ if $a \leq b$.
From now on, 
we set
\[m=n(L-1)
\] 
as well as in Remark~\ref{remark:bT}.

First
we shall construct the formulae for  $P_0^{(j)}$ $(0 \leq j \leq L-1)$.

\paragraph{}
\noindent
\underline{(i) Case $j=0$} \quad
The approximation condition (\ref{eq:spa}) requires that
\[
\left[
 f_i  P_0^{(0)}
\right]^{m+n-1}_{m}=0 \quad (1 \leq i \leq L-1)
\]
since $P_i^{(0)}$ $(i \neq 0)$ is a polynomial of degree at most $m-1$.
If we write
\[
P_0^{(0)}(w)=\sum_{k=0}^{ m  } b_{k}w^k,
\]
then we find a system
\begin{equation}  \label{eq:lspa}
\begin{bmatrix}
A^1_{m}(n,m+1) 
\\
A^2_{m}(n,m+1) 
\\
\vdots
\\
A^{L-1}_{m}(n,m+1) 
\end{bmatrix}
\begin{bmatrix} 
b_0 \\b_1 \\ \vdots \\b_{m}
\end{bmatrix}
=
\left.
\begin{bmatrix} 
0  \\ \vdots \\0
\end{bmatrix}
\right\} \text{\footnotesize $m$}
\end{equation}
of homogeneous linear equations for the $m+1$ unknowns 
$b_0, \ldots, b_{m}$. 
The solution of (\ref{eq:lspa}) is unique up to multiplication by constants if and only if
the rank of the $m \times (m+1)$ matrix in the left-hand side equals $m$.

\paragraph{}
\noindent
\underline{(ii) Case $1 \leq j \leq L-1$} \quad
Similarly, it follows from (\ref{eq:spa}) that
\begin{align*}
\left[
 f_i  P_0^{(j)}
\right]^{m+n-1}_{m}&=0 \quad (1 \leq i\leq j-1)
\\
\left[
 f_j  P_0^{(j)}
\right]^{m+n-2}_{m}&=0 \quad  (i=j)
\\
\left[
 f_i  P_0^{(j)}
\right]^{m+n-2}_{m-1}&=0 \quad  (j+1 \leq i \leq L-1).
\end{align*}
These amount to the simultaneous linear equation
\[
\begin{bmatrix}
A^1_{m}(n,m) 
\\
\vdots
\\
A^{j-1}_{m}(n,m) 
\\
A^{j}_{m}(n-1,m) 
\\
A^{j+1}_{m-1}(n,m)
\\
\vdots
\\
A^{L-1}_{m-1}(n,m) 
\end{bmatrix}
\begin{bmatrix} 
b_0 \\b_1 \\ \vdots \\b_{m-1}
\end{bmatrix}
=
\left.
\begin{bmatrix} 
0 \\ \vdots \\0
\end{bmatrix}
\right\} \text{\footnotesize $m-1$}
\]
for the $m$ unknown coefficients $b_0, \ldots, b_{m-1}$ of the polynomial
\[
P_0^{(j)}(w)=\sum_{k=0}^{ m-1  } b_{k}w^k.
\]

\begin{prop} 
\label{prop:detp}

The polynomials $P^{(j)}_0$ admit
the following 
determinantal representations{\rm:}
\[
P_0^{(0)}(w) 
=\frac{1}{{\rm NP}^{(0)}}
\det 
\begin{bmatrix}
1,w,w^2, \ldots, w^{m}
\\
A^1_{m}(n,m+1) 
\\
\vdots
\\
A^{L-1}_{m}(n,m+1) 
\end{bmatrix}
\]
and 
\[
P_0^{(j)}(w) 
=\frac{1}{{\rm NP}^{(j)}}
\det 
\begin{bmatrix}
1,w,w^2, \ldots, w^{m-1}
\\
A^1_{m}(n,m) 
\\
\vdots
\\
A^{j-1}_{m}(n,m)
\\ 
A^{j}_{m}(n-1,m)
\\ 
A^{j+1}_{m-1}(n,m)
\\
\vdots \\
A^{L-1}_{m-1}(n,m) 
\end{bmatrix}
\]
for $1 \leq j \leq L-1$,
where ${\rm NP}^{(j)}$ are some normalizing constants. 
\end{prop}

Next the other $P_i^{(j)}$ ($i\neq 0$) can be written as follows:
\[
\begin{array}{lllll}
\bullet \ \text{if $j=0$} &&
P_i^{(0)} 
= 
\left[
 f_i  P_0^{(0)}
\right]^{m-1}_{0};
\\ 
\bullet \ \text{if $1 \leq  j \leq L-1$} &&
P_i^{(j)} = 
\left[
 f_i  P_0^{(j)}
\right]^{m-1}_{0}  
&\text{for $1 \leq i \leq j-1$},
\\
&&P_j^{(j)} = 
w \left[
 f_j  P_0^{(j)}
\right]^{m-1}_{0}
&\text{for  $i =j$},
\\
&&P_i^{(j)} = 
w \left[
 f_i  P_0^{(j)}
\right]^{m-2}_{0}
&\text{for  $j+1 \leq i \leq L-1$}.
\end{array}
\]
We will choose the normalization so that 
the diagonal part
$P^{(i)}_i(w)$ becomes monic.
Accordingly, we obtain
\begin{equation}
\label{eq:np0}
{\rm NP}^{(0)}
=\det 
\begin{bmatrix}
0, \ldots,0,1
\\
A^1_{m}(n,m+1) 
\\
\vdots
\\
A^{L-1}_{m}(n,m+1) 
\end{bmatrix}
= 
(-1)^m\det 
\begin{bmatrix}
A^1_{m}(n,m) 
\\
\vdots
\\
A^{L-1}_{m}(n,m) 
\end{bmatrix}
=(-1)^{\frac{m(m-n)}{2} +n(L-1)} \Delta^{(L)}
\end{equation}
and
\begin{equation}
\label{eq:npk}
{\rm NP}^{(j)}
=
\det 
\begin{bmatrix}
a^j_{m-1},a^j_{m-2}, \ldots, a^j_{0}
\\
A^1_{m}(n,m) 
\\
\vdots
\\
A^{j-1}_{m}(n,m)
\\ 
A^{j}_{m}(n-1,m)
\\ 
A^{j+1}_{m-1}(n,m)
\\
\vdots \\
A^{L-1}_{m-1}(n,m) 
\end{bmatrix}
=
(-1)^{n(j-1)}
\det 
\begin{bmatrix}
A^1_{m}(n,m) 
\\
\vdots
\\
A^{j-1}_{m}(n,m)
\\ 
A^{j}_{m-1}(n,m)
\\ 
\vdots \\
A^{L-1}_{m-1}(n,m) 
\end{bmatrix}
=(-1)^{\frac{m(m-n)}{2} +n(j-1)} \Delta^{(j)}
\end{equation}
for $1 \leq j \leq L-1 $.
Here we have employed
 the notation of the block-Toeplitz determinant
(see Remark~\ref{remark:bT}) in view of
(\ref{eq:tenchi}).

\begin{remark}\rm
If $i \geq j \geq 1$, 
$P^{(j)}_i(w)$ is divisible by $w$.
By construction (see Corollary~\ref{cor:R}), we observe that
$R^{-1}$ takes the form
(cf. (\ref{eq:R}))
\[
R^{-1}
= 
w^{-m}
\left(
\begin{bmatrix}
* &     & \\ 
\vdots  &$\mbox{\Huge 0}$ &  \\ 
*&&
\end{bmatrix} 
+O(w)
\right).
\]
\end{remark}

\section{From vector continued fractions to Schlesinger transformations}
\label{sect:cf}

In this section 
we present an alternative construction
 of the same Schlesinger transformation (considered in Sect.~\ref{sect:hpa})
through an algorithm for expanding a vector-valued function
into a 
vector continued fraction.

\subsection{Algorithm for vector continued fraction expansion}

Let 
${\boldsymbol f}={}^{\rm T}(f_0,\ldots,f_{L-1})$ 
be the $L$-tuple of formal power series (\ref{eq:fi}).
For simplicity, 
we assume tentatively  that $f_i(0) \neq 0$ for all $ i$.
We abbreviate the constant term $a^i_0$
of $f_i(w)=\sum_{j=0}^\infty a^i_j w^j$ as $a^i$.

First we apply a left multiplication of a permutation matrix
to ${\boldsymbol f}$:
\[
\begin{bmatrix}
0 &  &&1
\\
1&0&&
\\
 &\ddots& \ddots\\ 
& & 1&0
\end{bmatrix}
\begin{bmatrix}
f_0 \\ f_1 \\
\vdots    \\ 
f_{L-1}
\end{bmatrix}
=\begin{bmatrix}
f_{L-1} \\ f_0 \\
\vdots    \\ 
f_{L-2}
\end{bmatrix}.
\]
Next we eliminate the constant term of $f_i$ by a subtraction of 
constant multiple of $f_{i+1}$ for each $0 \leq i \leq L-2$
and by a multiplication by $w$ for $i=L-1$:
\[
 \begin{bmatrix}
w &  &&&&
\\
&1&-\frac{a^0}{a^1}&&&
\\
&&1&-\frac{a^1}{a^2}&&
\\
&& &\ddots& \ddots
\\ 
&&&&1&-\frac{a^{L-3}}{a^{L-2}}
\\
-\frac{a^{L-2}}{a^{L-1}}& & &&&1
\end{bmatrix}
\begin{bmatrix}
f_{L-1} \\ f_0 \\ f_1 \\
\vdots    \\ f_{L-3} \\
f_{L-2}
\end{bmatrix}
=
O(w)
=
w
\begin{bmatrix}
f_0' \\  f_1' \\
\vdots    \\ 
 f_{L-1}'
\end{bmatrix}.
\] 
Eventually
we obtain a new $L$-tuple
 ${\boldsymbol f'}={}^{\rm T}(f_0',\ldots,f_{L-1}')$ 
of power series from the original ${\boldsymbol f}$.

The above procedure
 is summarized as a left multiplication
\begin{equation} \label{eq:f'}
{\boldsymbol f'}=T
{\boldsymbol f}
\end{equation}
of an invertible matrix
\[
T=
\frac{1}{w}
 \begin{bmatrix}
0 &  &&&&w
\\
1&-\frac{a^0}{a^1}&&&&
\\
&1&-\frac{a^1}{a^2}&&&
\\
& &\ddots& \ddots&
\\ 
&&&1&-\frac{a^{L-3}}{a^{L-2}}&
\\
 & &&&1&-\frac{a^{L-2}}{a^{L-1}}
\end{bmatrix}
\]
with $\det T=(-w)^{1-L}$. 
This is
an analogue of the Euclidean algorithm
and can be repeated generically. 
Let 
${\boldsymbol f}[k]={}^{\rm T}(f_0[k],\ldots,f_{L-1}[k])$  
denote the corresponding vector of power series at the $k$th step and let $a^i[k]$ 
denote their constant terms.
Hence, we have
\begin{equation}  \label{eq:fn}
{\boldsymbol f}[k+1]=T[k]{\boldsymbol f}[k]
\quad \text{and}
\quad
{\boldsymbol f}[0]={\boldsymbol f},
\end{equation}
where
\begin{equation}
\label{eq:Tn}
T[k]=\frac{1}{w}
 \begin{bmatrix}
0 &  &&&&w
\\
1&-\frac{a^0[k]}{a^1[k]}&&&&
\\
&1&-\frac{a^1[k]}{a^2[k]}&&&
\\
& &\ddots& \ddots&
\\ 
&&&1&-\frac{a^{L-3}[k]}{a^{L-2}[k]}&
\\
 & &&&1&-\frac{a^{L-2}[k]}{a^{L-1}[k]}
\end{bmatrix}.
\end{equation}

On the other hand,
solving (\ref{eq:f'})
 for $\boldsymbol f$ yields
\[
f_i= w \sum_{j=i}^{L-1}
\frac{a^i}{a^j}f'_{j+1}
\quad (0 \leq i \leq L-1)
\]
where $f'_L=f'_0/w$.
Let us introduce the 
inhomogeneous coordinates 
${\boldsymbol \varphi}={}^{\rm T}(\varphi_1,\ldots,\varphi_{L-1})$ by $\varphi_i=f_i/f_0$.
Therefore, we have
\begin{equation}  \label{eq:phi}
\frac{\varphi_i}{\varphi_{L-1}}=
\frac{a^i}{a^{L-1}}+
w \sum_{j=i}^{L-2} \frac{a^i}{a^j} \varphi'_{j+1}
\quad (0 \leq i  \leq L-2)
\end{equation}
where $\varphi_0=1$.

\begin{dfn}[cf. \cite{ns91, par81}]\rm
Let
${\boldsymbol \varphi}={}^{\rm T}(\varphi_1,\ldots,\varphi_{L-1})$ 
be an $(L-1)$-vector
such that $\varphi_1\neq 0$.
Then the vector
\[
\iota({\boldsymbol \varphi})=\frac{1}{{\boldsymbol \varphi}}=
{}^{\rm T}
\left( \frac{\varphi_2}{\varphi_1}, \frac{\varphi_3}{\varphi_1}, \ldots,  \frac{\varphi_{L-1}}{\varphi_1},  \frac{1}{\varphi_1} \right)
\]
is called the {\it reciprocal of ${\boldsymbol \varphi}$}.
Note that  $\iota^L=\id$.
\end{dfn}

Under this notation, 
the correspondence (\ref{eq:phi}) can be translated into
\[
{}^{\rm T}\left(\frac{1}{\varphi_{L-1}}, \frac{\varphi_1}{\varphi_{L-1}}, \ldots,\frac{\varphi_{L-2}}{\varphi_{L-1}}\right)
=\iota^{-1}( {\boldsymbol \varphi} )=
{\boldsymbol a}+ w B {\boldsymbol \varphi'}
\]
or equivalently into
\[
{\boldsymbol \varphi}=\frac{1}{{\boldsymbol a}+ w B{\boldsymbol \varphi'} },
\]
where
\[
{\boldsymbol a}= 
\frac{1}{a^{L-1}}
\begin{bmatrix}
a^0 \\ a^1\\
\vdots    \\ 
a^{L-2}
\end{bmatrix}
\quad \text{and} \quad
B= \begin{bmatrix}
1 & \frac{a^0}{a^1} & \frac{a^0}{a^2} &\cdots & \frac{a^0}{a^{L-2}}  
\\ 
& 1 & \frac{a^1}{a^2} &\cdots & \frac{a^1}{a^{L-2}}  
\\
&&\ddots&\ddots&\vdots
 \\ 
&&& 1 & \frac{a^{L-3}}{a^{L-2}}
\\
&&&&1
\end{bmatrix}.
\]
Namely,
the vector
${\boldsymbol a}$ is determined as
the constant term of $\iota^{-1}({\boldsymbol \varphi})$
and the matrix $B$ is then specified by ${\boldsymbol a}$.
Let 
${\boldsymbol \varphi}[k]$ denote the inhomogeneous coordinates 
of the vector ${\boldsymbol f}[k] \in  {\mathbb C}^{L}[\![w]\!]$
at the $k$th step.
Taking the reciprocal repeatedly in this way,
we obtain formally the {\it vector continued fraction}
\begin{equation} \label{eq:stiel}
{\boldsymbol \varphi}= \cfrac{1}{{\boldsymbol a}[0] + \cfrac{ w B[0]}{ {\boldsymbol a}[1] + \cfrac{ w B[1]}{ {\boldsymbol a}[2] + \cfrac{w B[2]}{ {\boldsymbol a}[3] +\ddots} }  }},
\end{equation}
which is regarded as an $(L-1)$-dimensional generalization 
of the {\it Stieltjes-type} continued fraction.
Refer to \cite[Appendix A]{jt09} for a classification of continued fractions.
Our algorithm
differs from the other
known examples
such as 
the
Jacobi--Perron algorithm;
cf. \cite{par81, per07}. 
Note also that some dynamical system, like the Toda lattice, 
has been studied
based on the connection among 
the Jacobi--Perron algorithm, 
rational approximations and bi-orthogonal polynomials; see \cite{ka84}.

The following theorem can be verified straightforwardly through the above algorithm,
as well as the case of a Stieltjes-type continued fraction (i.e., $L=2$ case).

\begin{thm} 
The $k$th convergents
{\rm(}rational functions{\rm)}
\begin{align*}
{\boldsymbol \Pi}_1&=\cfrac{1}{{\boldsymbol a}[0]},
\\
{\boldsymbol \Pi}_k&= \cfrac{1}{{\boldsymbol a}[0] + \cfrac{ w B[0]}{ {\boldsymbol a}[1] + \cfrac{ w B[1]}{ {\boldsymbol a}[2] + \ddots + \cfrac{w B[k-2] }{ {\boldsymbol a}[k-1]  } }  }}
\in {\mathbb C}^{L-1}(w)
\quad (k \geq 2)
\end{align*}
of the vector continued fraction 
{\rm(\ref{eq:stiel})}
provide approximants of
the vector   
${\boldsymbol \varphi} \in {\mathbb C}^{L-1}[\![w]\!]$
of power series
in the sense that
${\boldsymbol \varphi}- {\boldsymbol \Pi}_k =O(w^k)$.
\end{thm}

In calculating ${\boldsymbol \Pi}_k$,
it is convenient to apply
the projective transformations
(\ref{eq:Tn})
 successively
as follows:
\[
\begin{bmatrix}
\varpi_0 \\ \varpi_1 \\
\vdots \\ \varpi_{L-1}
\end{bmatrix}
=
T[0]^{-1}T[1]^{-1} \cdots T[k-1]^{-1} 
\begin{bmatrix}1 \\0 \\
\vdots \\ 0
\end{bmatrix}
\quad \text{and} \quad
{\boldsymbol \Pi}_k=\frac{1}{\varpi_{0}}\begin{bmatrix}\varpi_1 \\ \varpi_2 \\
\vdots \\ \varpi_{L-1}
\end{bmatrix}.
\]

\subsection{Schlesinger transformations, revisited}

Let $Y=Y(z)$ be a solution 
(\ref{eq:fss}) 
of the Fuchsian system
(\ref{eq:fs})
having the local behaviors
\begin{align*}
Y &=
\Psi (z) \cdot {\rm diag \ }  (z^{\varepsilon_{0,j}})_{0 \leq j \leq L-1} 
\qquad \text{(near $z=0$)}
\\
&=
\Phi (w) \cdot {\rm diag \ }  (w^{\varepsilon_{\infty,j}})_{0 \leq j \leq L-1}
 \cdot C
 \qquad \text{(near $z=1/w=\infty$)}
\end{align*}
where the power series parts are normalized by
\[
\Psi(z)
=
\begin{bmatrix}
1 &  \cdots &*\\  
&\ddots&  \vdots \\
&&1
\end{bmatrix} 
+O(z),
\quad
\Phi(w)
=
\begin{bmatrix}
* &  &\\  
 \vdots&\ddots& \\ 
*& \cdots&*
\end{bmatrix}
+O(w)
\]
and $C$ is the connection matrix.
Let
${\boldsymbol f}={}^{\rm T}(f_0, \ldots,f_{L-1})$
denote the first column of $\Phi(w)$.

It is clear from the construction of the matrix $T=T[0]$
that
\[
T \Phi=
\left(\begin{bmatrix}
* &  &\\  
 \vdots&\ddots& \\ 
*& \cdots&*
\end{bmatrix}
+O(w)
\right) \cdot 
 {\rm diag \ } 
(1, w^{-1}, \ldots,w^{-1})
\]
near $w=0$ ($z=\infty$).
On the other hand, it holds that
\begin{align*}
T \Psi
&=
z
\left(
\begin{bmatrix}
*&*&\cdots&*&1 \\
1 & *&\cdots &*&0\\  
&1&\ddots& \vdots&\vdots  \\
&&\ddots& * & 0 \\
&&&1&0
\end{bmatrix} 
+O(z)
\right)
 \cdot 
 {\rm diag \ } 
 (1, \ldots,1, z^{-1})
\\
&=
z
\left(\begin{bmatrix}
1 & \cdots &*\\  
&\ddots& \vdots \\
&&1
\end{bmatrix} 
+O(z)
\right) 
\begin{bmatrix}
0&&& z^{-1}
\\
1 &  0&\\  
&\ddots& \ddots \\
&&1&0
\end{bmatrix} 
\end{align*}
near $z=0$.
After we repeat the same procedure $L$ times,
the power series part thus recovers its original form:
\[
T[L-1]\cdots T[1]T[0]\Psi
=
z^{L-1}
\left(\begin{bmatrix}
1 & \cdots &*\\  
&\ddots& \vdots \\
&&1
\end{bmatrix} 
+O(z)
\right). 
\] 
In conclusion, 
the matrix 
$z^{1-L} T[L-1]\cdots T[1]T[0]$
turns out to be a polynomial in $z$ and to be
the multiplier of the Schlesinger transformation
shifting the characteristic exponents at $z=\infty$ by
$(L-1,-1,\ldots,-1)$;
cf. Theorem~\ref{thm:schl}.

\begin{remark}\rm
In fact,
the same approximation problem considered in Sect.~\ref{sect:hpa} 
appears in the following manner.
Concerning the polynomial matrix $w^LT[L-1]\cdots T[1]T[0] $,
we observe from the form of $T[k]$ that
the diagonal entries are all monic linear functions,
the strictly upper triangular part is linear and divisible by $w$,
and 
the strictly lower triangular part is a constant.  
Moreover,
in view of (\ref{eq:fn}) we have
\[
w^LT[L-1]\cdots T[1]T[0] {\boldsymbol f}=O(w^L),
\]
which coincides with the approximation condition (\ref{eq:hpa}), and thus
\[
w^LT[L-1]\cdots T[1]T[0] =
\begin{bmatrix} 
\widetilde{\boldsymbol Q}^{(0)}(w)
\\
\vdots
\\
\widetilde{\boldsymbol Q}^{(L-1)}(w)
\end{bmatrix}
\]
under $n=1$.
\end{remark}

\section{Application to  isomonodromic deformations}
\label{sect:aid}

In this section we first 
review some basic results on
 the {\it Schlesinger system},
which 
governs
isomonodromic deformations of a
Fuchsian system
of linear ordinary differential equations. 
As explained in Sect.~\ref{subsect:st},
a Schlesinger transformation preserves the monodromy of the 
Fuchsian system under consideration 
and, thereby, 
leads to a discrete symmetry of the associated Schlesinger system.
Combining this fact with the result in Sect.~\ref{sect:det}
reveals a determinantal nature of
isomonodromic deformations.
Next we 
treat
a particular case of 
the Schlesinger systems
unifying various Painlev\'e-type differential 
equations
and show its relationship with certain hypergeometric functions,
which will be needed later.

\subsection{Schlesinger systems and their symmetries}
\label{subsect:ss}

Let us consider again the $L \times L$ Fuchsian system (\ref{eq:fs}):
\begin{equation}
\label{eq:fs2}
\frac{{\rm d} Y}{{\rm d} z}
=A Y
=\sum_{i=0}^{N+1} \frac{A_i}{z-u_i} Y,
\end{equation}
where $u_0=1$ and $u_{N+1}=0$.
We start with 
a well-known result on isomonodromic deformations
of (\ref{eq:fs2}).

\begin{thm} 
The monodromy of a fundamental solution $Y$,
i.e. $\det Y \neq 0$,
does not depend on ${\boldsymbol u}=(u_1,\ldots,u_N)$
if and only if 
\begin{equation} \label{eq:defe}
B_i=\frac{\partial Y}{\partial u_i} Y^{-1}
\quad
(1 \leq i \leq N)
\end{equation}
are rational functions in $z$.
\end{thm}

We henceforth impose
on our Fuchsian system (\ref{eq:fs2}) 
 the following assumptions: 
\begin{quote}
(i) all the residue matrices $A_i$ are  {\it semi-simple}, i.e. diagonalizable;
\\
(ii) there is no integer difference other than zero
among the eigenvalues of each $A_i$.
\end{quote}
Let us choose a normalization as before 
such that
$A_{N+1}$ and $A_{N+2}=-\sum_{i=0}^{N+1} A_i$ are upper and lower triangular matrices, respectively.  
Then
we can take a fundamental solution 
$Y=Y(z)$
of the form 
\begin{align*}
Y &=
\Psi (z) \cdot {\rm diag \ }  (z^{\varepsilon_{0,j}})_{0 \leq j \leq L-1} 
\qquad \text{(near $z=0$)}
\\
&=
\Phi (w) \cdot {\rm diag \ }  (w^{\varepsilon_{\infty,j}})_{0 \leq j \leq L-1}
 \cdot C
 \qquad \text{(near $z=1/w=\infty$)}
\end{align*}
with
\begin{align*}
\Psi(z)
&=
\begin{bmatrix}
1 &  \cdots &*\\  
&\ddots& \vdots \\
&&1
\end{bmatrix} 
+O(z),
\\
\Phi(w)
&
=
\begin{bmatrix}
* &   &\\  
\vdots&\ddots&  \\
*& \cdots &*
\end{bmatrix} 
+O(w)
=
\left(
\begin{bmatrix}
1 &  &\\  
 \vdots&\ddots& \\ 
*& \cdots&1
\end{bmatrix}
+O(w)
\right)
\Xi
\end{align*}
and
$\Xi={\rm diag \ }  (\xi_0, \ldots, \xi_{L-1})$,
where each $\xi_j \neq 0$ may depend on 
$\boldsymbol u=(u_1, \ldots,u_N)$.
Moreover, 
we have
\begin{equation}
Y=G_i
(I+O(z-u_i))
(z-u_i)^{\Lambda_i}
\label{eq:solatui}
\end{equation}
near each of the other regular
singularities $z=u_i$ $(0 \leq i \leq N)$,
where
$G_i$ and $\Lambda_i$ are certain constant matrices
satisfying
$G_i \Lambda_i {G_i}^{-1}=A_i$.
Therefore, 
the monodromy matrices of $Y$ attached to loops around 
$z=u_i$ $(0 \leq i \leq N)$ and $z=0, \infty$
read
\[
e^{2 \pi \sqrt{-1}  \Lambda_i}, \quad 
e^{2 \pi \sqrt{-1} 
 {\rm \ diag \ }  (\varepsilon_{0,j})_{0 \leq j \leq L-1} },
 \quad
C^{-1}
e^{2 \pi \sqrt{-1} 
 {\rm \ diag \ }  (\varepsilon_{\infty,j})_{0 \leq j \leq L-1} }
C.
\] 

Suppose now that every monodromy matrix
of $Y$ 
is constant with respect to $\boldsymbol u$ and, additionally, 
so is the connection matrix $C$. 
Then
the rational functions $B_i=B_i(z)$
can be explicitly written as 
\begin{equation} \label{eq:bi}
B_i=\frac{A_i}{u_i-z}-\frac{1}{u_i}
(A_i)_{\rm LT},
\end{equation}
where
$(A_i)_{\rm LT}$ denotes the lower triangular part
of $A_i$;
see Appendix for details.
Note in particular that the diagonal part 
$(A_i)_{\rm D}$ of $A_i$ 
is expressible in terms of 
$\Xi$ as
\begin{equation} \label{eq:Aid}
(A_i)_{\rm D}= - u_i \frac{\partial}{\partial u_i}
\log \Xi.
\end{equation}
The compatibility condition
\[\frac{\partial A}{\partial u_i} -\frac{\partial B_i}{\partial z}+[A,B_i]=0
\]
 of  
 (\ref{eq:fs2}) and (\ref{eq:defe})
 is equivalent to a set of nonlinear differential equations for 
 the matrices
 $A_i$
with respect to ${\boldsymbol u}$,
which is called  
the {\it Schlesinger system}
\cite{sch12}.
If $(L,N)=(2,1)$,
then the Schlesinger system reduces to the sixth Painlev\'e equation
$P_{\rm VI}$.

Next we shall investigate how 
the solution $Y$ and the coefficient $A=A(z)$ of the Fuchsian system (\ref{eq:fs2})
are connected with each other.
Concerning the power series expansion
\[
Y=
\Phi (w) \cdot {\rm diag \ }  (w^{\varepsilon_{\infty,j}}) \cdot C,
\quad \Phi(w)= \sum_{k=0}^\infty \Phi_k w^k
 \]
 at the point of infinity ($z=1/w=\infty$),
 the coefficients $\Phi_k$ turn out to be
 polynomials in the off-diagonal entries of 
 the lower triangular matrix 
 $A_{N+2}=-\sum_{i=0}^{N+1}A_i$
through Frobenius' method.
Conversely, substituting this solution $Y$ in (\ref{eq:fs2}),
we find that 
\begin{align*}
z \frac{{\rm d}Y}{{\rm d} z}
&=-w \frac{{\rm d}Y}{{\rm d} w}
\\
&=-w \frac{{\rm d}}{{\rm d} w} 
\left(\Phi_0+\Phi_1 w+\Phi_2 w^2 +\cdots
\right)
{\rm diag \ }
(w^{\varepsilon_{\infty,j}})
\cdot C
\\
&=
-
\left[
\left(\Phi_0+\Phi_1 w+\Phi_2 w^2 +\cdots \right)
{\rm diag \ }
(\varepsilon_{\infty,j})
\right.
\\
&\qquad \qquad
\left.
+ \Phi_1 w +2 \Phi_2 w^2 +3 \Phi_3 w^3 +\cdots  \right]
{\rm diag \ }
(w^{\varepsilon_{\infty,j}})
\cdot C
\end{align*}
and
\begin{align*}
z \sum_{i=0}^{N+1}
\frac{A_i}{z-u_i} Y 
&=\sum_{i=0}^{N+1}
\frac{A_i}{1-u_iw} Y 
\\
&=\left[
A_{N+1}+
\sum_{i=0}^{N} \left(1+ u_iw 
+{u_i}^2 w^2  + \cdots\right)A_i
\right]Y.
\end{align*}
Recall  $u_{N+1}=0$ here.
Equating the coefficients of $w^k$ in these power series yields
\[
A_{N+2}
=\Phi_0  \cdot 
 {\rm diag \ }  (\varepsilon_{\infty,j})
 \cdot {\Phi_0}^{-1}
 \quad \text{for $k=0$,}
\]
thus $A_{N+2}$ becomes a polynomial in the entries of the leading coefficient $\Phi_0$ of the solution,
and also
\[
\sum_{i=0}^{N} u_i A_i
=\left(
-\Phi_1  \cdot {\rm diag \ } (\varepsilon_{\infty,j}+1) +
\Phi_0 \cdot {\rm diag \ }(\varepsilon_{\infty,j}) \cdot {\Phi_0}^{-1} \Phi_1 
\right){\Phi_0}^{-1}
\quad \text{for $k=1$.}
\]
Moreover, if one needs  
similar expressions for all other residue matrices $A_i$
besides $A_{N+2}$,
it is convenient to use the {\it deformation equation}
(\ref{eq:defe});
one can verify in fact
\[
A_i 
=\left(
\frac{\partial \Phi_0}{\partial u_i}
{\Phi_0}^{-1}\Phi_1-
\frac{\partial \Phi_1}{\partial u_i}
\right)
{\Phi_0}^{-1}\quad (1 \leq i \leq  N).
\]
In summary, 
each residue matrix $A_i$ of $A(z)$ is expressible as a polynomial in
the entries of the coefficients of $Y$ (and their derivatives with respect to ${\boldsymbol u}$),
and vice versa.

A Schlesinger transformation 
keeps the monodromy of the Fuchsian system invariant
but shifts its characteristic exponents by integers;
recall Sect.~\ref{subsect:st}.
Consequently, it gives rise to a discrete symmetry of the Schlesinger system
via the above correspondence between 
the solutions and coefficients of the Fuchsian system.
On the other hand,
any ingredient of Schlesinger transformations or of 
the associated rational approximations 
is described in terms of block-Toeplitz determinants;
recall Sect.~\ref{sect:det}.
This fact thus provides a natural explanation for 
the determinantal structure appearing in solutions of 
isomonodromic deformations, e.g. 
Painlev\'e equations.
Refer to \cite{imt15} for a detailed investigation of the determinantal structure in Jimbo--Miwa--Ueno's {\it $\tau$-functions} \cite{jmu81}
for a general framework 
admitting irregular singularities.

\subsection{Polynomial Hamiltonian system ${\cal H}_{L,N}$ of isomonodromy type}

We turn now to a particular case of the Schlesinger systems,
which will be the main subject in the rest of this paper.

Consider an $L \times L$ Fuchsian system 
of the form (\ref{eq:fs2})
whose {\it spectral type} 
is given by the 
partitions of $L$:
\[
\begin{array}{ll}
1, L-1
&
\text{at $z=u_i$ $(0 \leq i \leq N )$ and} 
\\
\underbrace{1,1, \ldots,1}_{L}
&
\text{at $z=0, \infty$,}
\end{array}
\]
which indicate how the characteristic exponents overlap
at each of the $N+3$ singularities.
Fix the characteristic exponents as listed in the following table (Riemann scheme):
\begin{equation} 
\label{eq:rs}
\begin{array}{|c|c|}
\hline
\text{Singularity} & \text{Characteristic exponents}  \\ \hline
u_i  \mbox{\ } (0 \leq i \leq N )& (-\theta_i, 0, \ldots,0)  \\ \hline
u_{N+1}=0 & (e_0,e_1,\ldots,e_{L-1})  \\ \hline    
u_{N+2}=\infty & (\kappa_{0}-e_0,\kappa_{1}-e_1, \ldots,\kappa_{L-1}-e_{L-1} )  \\ 
\hline
   \end{array}
\end{equation}
Assume the sum of all the characteristic exponents equals zero
(Fuchs relation), 
i.e.
\begin{equation}
\label{eq:norm1}
\sum_{l=0}^{L-1}\kappa_l
=\sum_{i=0}^N \theta_i.
\end{equation}
Let
$A_{N+1}$ and $A_{N+2}=-\sum_{i=0}^{N+1} A_i$
be upper and lower triangular matrices, respectively.
Then such a Fuchsian system,
denoted by  ${\cal L}_{L,N}$,
can be
 parametrized as follows:
 \begin{align}
A_i 
&=  {}^{\rm T}\left(b_0^{(i)}, b_1^{(i)}, \ldots ,b_{L-1}^{(i)}\right)
\cdot
\left(c_0^{(i)}, c_1^{(i)}, \ldots ,c_{L-1}^{(i)}\right)\quad \text{with} \quad 
c_0^{(i)}=1
\quad
(0 \leq i \leq N),
 \nonumber
\\
A_{N+1} 
&= 
 \begin{bmatrix}
   e_0 &   w_{0,1}      & \cdots  & w_{0,L-1}
   \\
        & e_1 &  \ddots       & \vdots
   \\     
          &        &\ddots& w_{L-2,L-1}
   \\
          &        &          & e_{L-1}    
\end{bmatrix},
 \label{eq:resmat}
\end{align}
under the relations
\[
({\rm tr \,} 
A_i =)-\theta_i=
  \sum_{k=0}^{L-1} b_k^{(i)}c_k^{(i)},
 \quad
 \kappa_l=
 - \sum_{i=0}^{N} b_l^{(i)} c_l^{(i)}
 \quad
\text{and} \quad 
w_{k,l}=-\sum_{i=0}^N b_k^{(i)} c_l^{(i)}
\quad (k < l);
\]
the last two of which
come from 
the triangularity of $A_{N+2}$.
Also,
we can and will normalize the characteristic exponents
at $z=0$
by
\begin{equation}
\label{eq:norm2}
{\rm tr \,} 
A_{N+1}=\sum_{k=0}^{L-1} e_k=\frac{L-1}{2}
\end{equation}
without loss of generality.

As shown in \cite{tsu14a}, 
the Schlesinger system 
governing isomonodromic deformations of ${\cal L}_{L,N}$
reduces to 
the multi-time Hamiltonian system
${\cal H}_{L,N}$:
\[
\frac{\partial q_k^{(i)}}{\partial x_j}=\frac{\partial H_j}{\partial p_k^{(i)}},
\quad 
\frac{\partial p_k^{(i)}}{\partial x_j}=-\frac{\partial H_j}{\partial q_k^{(i)}}
\quad 
\left(
\begin{array}{c}
 1 \leq i , j \leq N \\
1 \leq k \leq L-1
\end{array}
\right).
\]
Here we let $x_i=1/u_i$ and define 
the Hamiltonian function $H_i$ by
\[
x_i H_i=
\sum_{k=0}^{L-1} e_k q_k^{(i)} p_k^{(i)} 
+\sum_{j=0}^N \sum_{0 \leq k < l \leq L-1}
 q_k^{(i)}   p_k^{(j)} q_l^{(j)} p_l^{(i)} 
+\sum^N_{
\begin{subarray}{c} 
j=0  \\
j \neq i
\end{subarray}
}
\frac{x_j}{x_i-x_j} 
 \sum_{k,l=0}^{L-1}
 q_k^{(i)}   p_k^{(j)} q_l^{(j)} p_l^{(i)}
\]
with
$x_0=q_l^{(0)}
= q_0^{(i)}=1$, 
$p_l^{(0)}
=
\kappa_{l} - \sum_{i=1}^{N} q_l^{(i)} p_l^{(i)}$
and 
$p_0^{(i)}
=\theta_i-\sum_{k=1}^{L-1} q_k^{(i)} p_k^{(i)}$.
Therefore, 
$H_i$ is a polynomial in the 
unknowns ({\it canonical variables})
\begin{equation}
\label{eq:canvar}
q_k^{(i)}=\frac{c_k^{(i)}}{c_k^{(0)}} \quad \text{and} \quad 
p_k^{(i)}=-b_k^{(i)}c_k^{(0)}
\quad 
\left(
\begin{array}{c}
 1 \leq i \leq N \\
1 \leq k \leq L-1
\end{array}
\right).
\end{equation}
The number of 
the constant parameters
\begin{equation}
\label{eq:const}
({\boldsymbol e},{\boldsymbol \kappa},{\boldsymbol \theta})
=(e_0,\ldots,e_{L-1},\kappa_0,\ldots,\kappa_{L-1},\theta_0,\ldots,\theta_N)
\end{equation}
contained in ${\cal H}_{L,N}$
 is essentially $2L+N-1$
 in view of (\ref{eq:norm1}) and (\ref{eq:norm2}).
For example, the case where $L=2$ and any $N \geq 1$
coincides with the 
Garnier system in $N$-variables
and, thereby, the first nontrivial case ${\cal H}_{2,1}$ does with the Hamiltonian form of $P_{\rm VI}$.

\begin{remark}\rm
We have 
{\it a priori} known
from their spectral type 
that the Fuchsian systems equipped with
the Riemann scheme (\ref{eq:rs})
constitute a $2N(L-1)$-dimensional family.
The coordinates of such a family are called 
{\it accessory parameters},
which are realized by the $2N(L-1)$
canonical variables (\ref{eq:canvar})
in this instance.
\end{remark}

\subsection{Solution of ${\cal H}_{L,N}$ in terms of hypergeometric function $F_{L,N}$}

\label{subsect:hg}

Although the phase space of ${\cal H}_{L,N}$
is a quite-complicated algebraic variety
in general,
there exists
a family of
solutions
parametrized by 
a point in the projective space
${\mathbb P}^{N(L-1)}$
when the constants $({\boldsymbol e},{\boldsymbol \kappa},{\boldsymbol \theta})$
take certain special values.
In fact, these solutions 
are written in terms of the {\it hypergeometric function} 
\begin{equation}
\label{eq:hg}
F_{L,N}
\left[
\begin{array}{c}
{\boldsymbol \alpha}, {\boldsymbol \beta} \\
{\boldsymbol  \gamma}
\end{array}
; {\boldsymbol x}
\right]
=
\sum_{m_i \geq0}
\frac{(\alpha_1)_{|{\boldsymbol m}|} \cdots (\alpha_{L-1})_{|{\boldsymbol m}|} (\beta_1)_{m_1} \cdots  (\beta_N)_{m_N} }{ (\gamma_1)_{| {\boldsymbol m}|} \cdots (\gamma_{L-1})_{| {\boldsymbol m}|} } 
\frac{{x_1}^{m_1} \cdots {x_N}^{m_N}}{{m_1}! \cdots {m_N}!},
\end{equation}
where
$|{\boldsymbol m}|=m_1+ \cdots +m_N$
and
$(a)_k=\Gamma(a+k)/\Gamma(a)$.
If $(L,N)=(2,1)$, then (\ref{eq:hg}) is exactly
Gau{\ss}'s 
hypergeometric function.

To state the result precisely, 
we introduce
the integral representation
\begin{equation}
\label{eq:euler}
y_0=\int_{c}
 \frac{U({\boldsymbol t}) \ 
 {\rm d}t_1 \cdots {\rm d}t_{L-1} }{ \prod_{l=1}^{L-1} (t_{l-1}-t_l)}
 = 
 \prod_{l=1}^{L-1} \frac{\Gamma(\alpha_l) \Gamma(\gamma_l-\alpha_l)}{\Gamma(\gamma_l)}
 \times F_{L,N}
\left[
\begin{array}{c}
{\boldsymbol \alpha}, {\boldsymbol \beta} \\
{\boldsymbol  \gamma}
\end{array}
; {\boldsymbol x}
\right]
\end{equation}
of $F_{L,N}$,
where 
the multi-valued function $U=U({\boldsymbol t})$ in ${\boldsymbol t}=(t_1,t_2,\ldots,t_{L-1})$
is given by
\[
U({\boldsymbol t})=
\prod_{l=1}^{L-1} {t_l}^{\alpha_l-\gamma_{l+1}}  
(t_{l-1}-t_l)^{\gamma_l-\alpha_l}
\prod_{i=1}^{N} (1-x_i t_{L-1})^{-\beta_i}
\quad \text{and} \quad t_0=\gamma_{L}=1,
\]
and the cycle $c$ is chosen to be
an $(L-1)$-simplex
\begin{equation}
\label{eq:simplex}
\{0 \leq t_{L-1} \leq \cdots \leq t_2 \leq t_1 \leq 1 \} \subset {\mathbb R}^{L-1}.
\end{equation}
Also, we introduce supplementarily
the integrals
\[
y_k^{(i)}
=\int_{c}
 \frac{U({\boldsymbol t})
 \ 
 {\rm d}t_1 \cdots {\rm d}t_{L-1} }{ 
 (x_i t_{L-1}-1)\prod_{\begin{subarray}{l} l=1 
 \\
  l \neq k \end{subarray}}^{L-1}(t_{l-1}-t_l)
 }
 \quad 
\left(
\begin{array}{c}
 1 \leq i  \leq N \\
1 \leq k \leq L-1
\end{array}
\right).
\]
We are now ready to state the hypergeometric solution of ${\cal H}_{L,N}$;
see \cite[Theorem~3.2]{tsu12b}.

\begin{thm} 
\label{thm:hgsol}
If $\kappa_0-\sum_{i=1}^N \theta_i=0$
then
the Hamiltonian system
${\cal H}_{L,N}$ 
possesses 
a solution 
\[
q_k^{(i)}=0 \quad \text{and} \quad
p_k^{(i)}
=
\theta_i\frac{y_k^{(i)}}{y_0}
\]
under the correspondence
\[
\alpha_k=e_k-e_0, \quad 
\beta_i=  - \theta_i, \quad  
\gamma_k= e_k-e_0-\kappa_k
\]
of constant parameters.
\end{thm}

The vector-valued function
\[
{\boldsymbol y}=
{}^{\rm T}\left(
y_0,y_1^{(1)},\ldots,y_{L-1}^{(1)},
y_1^{(2)},\ldots,y_{L-1}^{(2)},
 \ldots,
y_1^{(N)},\ldots,y_{L-1}^{(N)}
\right)
\]
satisfies a certain linear Pfaffian system
${\cal P}_{L,N}$
of rank $N(L-1)+1$,
whose
fundamental solution is prepared 
by collecting admissible cycles 
along with the foregoing $(L-1)$-simplex (\ref{eq:simplex}).
Note that the linear space of these cycles,
i.e. {\it twisted de Rham homology group}, 
is generated
by the chambers framed by the real section
of branch locus of 
$U=U({\boldsymbol t})$;
see  \cite{ak11, tsu12b}.
Of course, Theorem~\ref{thm:hgsol} is valid 
for any solution $\{y_0, y_k^{(i)}\}$
of ${\cal P}_{L,N}$.

The Fuchsian system ${\cal L}_{L,N}$
is specialized as
$\kappa_0-\sum_{i=1}^N \theta_i=0$ and
\begin{align}
&b_0^{(0)}=0,  \quad
b_k^{(0)}= \kappa_k y_0, \quad
 b_0^{(i)}=-\theta_i,
\quad
b_k^{(i)}= \theta_i y_k^{(i)}, 
\nonumber
\\
&c_k^{(0)}=\frac{-1}{y_0}, \quad
c_k^{(i)}=0 \quad  
\left(
\begin{array}{c}
 1 \leq i \leq N \\
1 \leq k \leq L-1
\end{array}
\right)
\label{eq:special}
\end{align}
along
the above hypergeometric solution of
${\cal H}_{L,N}$;
it thus becomes reducible.
In fact, via the gauge transformation
\[
Y= \left(w^{\kappa_0-e_0} \prod_{i=1}^N (1- u_i w)^{-\theta_i} {u_i}^{\theta_i} \right)Y',
\]
we have a solution of the form
\begin{align*}
Y' 
&=
\begin{bmatrix}
1 & 0 & \cdots & 0 \\
f_1 & & & \\
\vdots && \mbox{\huge $W$}  & \\
f_{L-1} &&&
\end{bmatrix}
\\
&=
\left(
\begin{bmatrix}
1   \\   
 \vdots&\ddots \\ 
*& \cdots&1
\end{bmatrix}
+O(w)
\right)
\cdot
{\rm diag \ }  (w^{\kappa_j -e_j -\kappa_0 + e_0 })_{0 \leq j \leq L-1},
\end{align*}
where
$f_k=f_k(w)$ $(1 \leq k \leq L-1)$
are holomorphic functions 
at $w=0$ defined by the integrals
\[
f_k
=
\int_{c}
 \frac{U({\boldsymbol t}) \ 
 {\rm d}t_1 \cdots {\rm d}t_{L-1} }{ 
 (1-w t_{L-1})
 \prod_{
 \begin{subarray}{l} 
 l=1 \\  l \neq k 
 \end{subarray}  
 }^{L-1}(t_{l-1}-t_l)
 }.
 \]
Note that an $(L-1) \times (L-1)$ matrix $W$  
can be described by Thomae's hypergeometric function 
$F_{L-1,1}$.
For details we refer to \cite{tsu14b},
in which 
a curious coincidence between 
${\cal P}_{L,N+1}$ and 
the {\it Lax pair} of ${\cal H}_{L,N}$, 
i.e.
the pair of
the original Fuchsian system 
(\ref{eq:fs2})
and its deformation equation
(\ref{eq:defe}),
is also discussed.

Our aim here is to generalize Theorem~\ref{thm:hgsol} by application of 
Schlesinger transformations starting from this hypergeometric solution
at $\kappa_0-\sum_{i=1}^N \theta_i=0$.
Notice that the Schlesinger transformation established in 
Theorem~\ref{thm:schl} 
shifts the constant parameters (\ref{eq:const})
as
\[
(\kappa_0, \kappa_1,\ldots, \kappa_{L-1})
\mapsto 
(\kappa_0+n(L-1), \kappa_1-n,\ldots, \kappa_{L-1}-n)
\]
while all the others are unchanged;
cf. the Riemann scheme (\ref{eq:rs}). 
Hence, by virtue of the algebraic relation
between the solution of ${\cal L}_{L,N}$
and the canonical variables of ${\cal H}_{L,N}$
(recall Sect.~\ref{subsect:ss} and also (\ref{eq:resmat}) and (\ref{eq:canvar})),
we know {\it in principle}
how to derive a solution 
$(\hat{q}_k^{(i)},\hat{p}_k^{(i)})$
of ${\cal H}_{L,N}$
at $\kappa_0-\sum_{i=1}^N \theta_i=n(L-1)$
for any positive integer $n$
even though the resulting expression in this way
will be terribly complicated.
In the next section
we explore this problem 
to achieve much simpler 
formulae for these special solutions.

\section{Solutions of  ${\cal H}_{L,N}$ in terms of
iterated hypergeometric integrals}
\label{sect:ihi}

This section is concerned with the Schlesinger transform of 
the hypergeometric solution of ${\cal H}_{L,N}$.
We present its explicit formula 
 by using the block-Toeplitz determinant whose entries are given by
the hypergeometric functions.
Key ingredients of the argument are 
the determinantal representations for the approximation polynomials;
see Sect.~\ref{sect:det}.
Moreover, we prove through Fubini's theorem and the Vandermonde determinant
that these block-Toeplitz determinants can be written
in the form of 
iterated hypergeometric integrals. 
Our result will be summarized in Theorem~\ref{thm:hgi},
which is regarded as a generalization of Theorem~\ref{thm:hgsol},
i.e. the previously known hypergeometric solution of ${\cal H}_{L,N}$.

\subsection{Preliminaries}
\label{subsect:prel}

Let 
$f_0(w), f_1(w), \ldots, f_{L-1}(w)$
be the functions defined by
\begin{align*}
f_0(w)&\equiv1,
\\
f_k(w)
&=
\int_{c}
 \frac{U({\boldsymbol t}) \ 
 {\rm d}t_1 \cdots {\rm d}t_{L-1} }{ 
 (1-w t_{L-1})
 \prod_{
 \begin{subarray}{l} 
 l=1 \\  l \neq k 
 \end{subarray}  
 }^{L-1}(t_{l-1}-t_l)
 }
  \quad (1 \leq k \leq L-1).
\end{align*}
If the cycle $c$ is chosen 
such that $|t_{L-1}|< \infty$,
then $f_k(w)$ 
is holomorphic at $w=0$.
For instance, it is 
enough to choose a bounded cycle as $c$. 
Accordingly, we have
a power series expansion
\[
f_k(w)= \sum_{j=0}^\infty h^k_j w^j
=h^k_0+h^k_1 w+h^k_2 w^2+\cdots
\]
with the coefficients
\begin{equation}
\label{eq:hjk}
h^k_j =
\int_{c}
 \frac{  {t_{L-1}}^{j} U({\boldsymbol t}) \ 
 {\rm d}t_1 \cdots {\rm d}t_{L-1} }{ 
 \prod_{
 \begin{subarray}{l} 
 l=1 \\  l \neq k 
 \end{subarray}  
 }^{L-1}(t_{l-1}-t_l)}
\end{equation}
for $1 \leq k \leq L-1$.
Observe that 
each $h^k_j $ can be 
regarded as a moment
\[
h^k_j =
\int_{ {\rm pr} (c) }
 {s}^{j} \ 
  {\rm d} \mu_k(s)
\]
of the `measure'
\[
 {\rm d} \mu_k(t_{L-1})
 =  
 \left(
 \int_{ c |t_{L-1} }
\frac{U({\boldsymbol t}) \ 
{\rm d}t_1 \cdots {\rm d}t_{L-2}}{
 \prod_{
 \begin{subarray}{l} 
 l=1 \\  l \neq k 
 \end{subarray}  
 }^{L-1}(t_{l-1}-t_l)  }
 \right)
{\rm d}t_{L-1}
\]
upon
 the following notations:
\begin{align*}
{\rm pr} (c)
&= \{ t_{L-1} \ | \ (t_1, \ldots, t_{L-2}, t_{L-1}) \in c
\},
\\
c |t_{L-1}
&= \{ 
(t_1, \ldots, t_{L-2})
 \ | \ (t_1, \ldots, t_{L-2}, t_{L-1}) \in c
\}.
\end{align*}
Namely
 $f_k(w)$  $(1 \leq k \leq L-1)$
 is  written as
the {\it Stieltjes transform}
\[
f_k(w)= 
\int_{{\rm pr} (c)}
 \frac{ {\rm d} \mu_k(s) }{ 
 1-w s}
\]
of a function $\mu_k=\mu_k(s)$.

In parallel, we introduce the functions
\begin{equation}
\label{eq:hj}
h_j^0 =
\int_{c}
 \frac{  {t_{L-1}}^{j} U({\boldsymbol t}) \ 
 {\rm d}t_1 \cdots {\rm d}t_{L-1} }{ 
 \prod_{ l=1 }^{L-1}(t_{l-1}-t_l)}
\end{equation}
also.
We thus see that
\[
h_j^0 =
\int_{ {\rm pr} (c) }
 {s}^{j} \ 
  {\rm d} \mu_0(s),
\]
where
\[
 {\rm d} \mu_0(t_{L-1})
 =  
 \left(
 \int_{ c |t_{L-1} }
\frac{U({\boldsymbol t}) \ 
{\rm d}t_1 \cdots {\rm d}t_{L-2}}{
 \prod_{ l=1 }^{L-1}(t_{l-1}-t_l)  }
 \right)
{\rm d}t_{L-1}.
\]

\begin{lemma}
The following linear relations {\rm(}contiguity relations{\rm)}
hold{\rm:} 
\begin{align}
\label{eq:cont1}
h^1_j+h^2_j+\cdots+h^{L-1}_j 
&=h^0_j-h^0_{j+1},
\\
\label{eq:cont2}
 h^k_{j}- x_i h^k_{j+1} 
 &= \ell_i (h^k_j),
\end{align}
where $\ell_i$ denotes the down-shift operator 
with respect to $\beta_i$
defined by 
\[
\ell_i(\beta_i)=\beta_i-1 \quad
\text{and} \quad
 \ell_i(\beta_j)=\beta_j \quad (i \neq j).
\]
\end{lemma}

\pf
By  
definitions (\ref{eq:hjk}) and (\ref{eq:hj})
it is immediate to verify these formulae. 
\qed
\\

\begin{remark} \rm
If $c$ is chosen to be the $(L-1)$-simplex
(\ref{eq:simplex}):
\[
\{0 \leq t_{L-1} \leq \cdots \leq t_2 \leq t_1 \leq 1 \} \subset {\mathbb R}^{L-1},
\]
both (\ref{eq:hjk}) and (\ref{eq:hj})
are
written by the hypergeometric function $F_{L,N}$.
Let us introduce the function
\[
h=h
\left[
\begin{array}{c}
{\boldsymbol \alpha}, {\boldsymbol \beta} \\
{\boldsymbol  \gamma}
\end{array}
; {\boldsymbol x}
\right]
=
\prod_{l=1}^{L-1} \frac{\Gamma(\alpha_l) \Gamma(\gamma_l-\alpha_l)}{\Gamma(\gamma_l)}
\times
F_{L,N}
\left[
\begin{array}{c}
{\boldsymbol \alpha}, {\boldsymbol \beta} \\
{\boldsymbol  \gamma}
\end{array}
; {\boldsymbol x}
\right].
\]
Then it holds by definition that
\[
h_j^k=
h\left[
\begin{array}{c}
\alpha_1+j+1, \ldots,\alpha_{k-1}+j+1, \alpha_{k}+j, \ldots, \alpha_{L-1}+j,
{\boldsymbol \beta}
\\
\gamma_1+j+1,\ldots,\gamma_{k}+j+1,\gamma_{k+1}+j, \ldots, \gamma_{L-1}+j
\end{array}
; {\boldsymbol x}
\right]
\quad (0 \leq k \leq L-1).
\]
For example we have
\begin{align*}
h^0_j &=
h\left[
\begin{array}{c}
\alpha_1+j, \ldots,\alpha_{L-1}+j,
{\boldsymbol \beta}
\\
\gamma_1+j,\ldots,\gamma_{L-1}+j
\end{array}
; {\boldsymbol x}
\right],
\\
h_j^1&=
h\left[
\begin{array}{c}
\alpha_{1}+j, \ldots, \alpha_{L-1}+j,
{\boldsymbol \beta}
\\
\gamma_1+j+1,\gamma_{2}+j, \ldots, \gamma_{L-1}+j
\end{array}
; {\boldsymbol x}
\right],
\\
h_j^2&=
h\left[
\begin{array}{c}
\alpha_{1}+j+1, \alpha_{2}+j,\ldots, \alpha_{L-1}+j,
{\boldsymbol \beta}
\\
\gamma_1+j+1,\gamma_{2}+j+1 ,\gamma_{3}+j, \ldots, \gamma_{L-1}+j
\end{array}
; {\boldsymbol x}
\right]
\end{align*}
and so forth.
\end{remark}

It is convenient to prepare the notation of the block-Toeplitz determinant
(cf. Remark~\ref{remark:bT})
for any $L$-tuple of nonnegative integers:
\[
{\boldsymbol n}=(n_0,n_1,\ldots,n_{L-1}) \in ({\mathbb Z}_{\geq 0})^L.
\]
Let $| {\boldsymbol n}|=\sum_{i=0}^{L-1} n_i$.
We set
\begin{align}
\Delta^{(k)}
({\boldsymbol n})
&:=
 \det
 \left[
\underbrace{
A^0_{n_0}(  | {\boldsymbol n} |,  n_0) \
\cdots  \
A^{k-1}_{n_{k-1}}(| {\boldsymbol n} |,  n_{k-1})}_{\text{$k$ blocks}}
 \
\underbrace{
A^{k}_{n_k-1}(| {\boldsymbol n} |,  n_k)  \
\cdots   \  
A^{L-1}_{n_{L-1}-1}(| {\boldsymbol n} |,  n_{L-1})
}_{\text{$L-k$ blocks}}
\right]
\nonumber
\\
&=
(-1)^{\sum_{i<j} n_in_j}
\det 
\begin{bmatrix}
A^0_{ | {\boldsymbol n} | }(n_0,| {\boldsymbol n} |) 
\\
\vdots
\\
A^{k-1}_{| {\boldsymbol n} |}(n_{k-1},| {\boldsymbol n} |)
\\ 
A^{k}_{| {\boldsymbol n} |-1}(n_k,| {\boldsymbol n} |)
\\ 
\vdots \\
A^{L-1}_{| {\boldsymbol n} |-1}(n_{L-1},| {\boldsymbol n} |) 
\end{bmatrix}
\label{eq:bT2}
\end{align}
for each $k$ $(0 \leq k \leq L)$.
If ${\boldsymbol n}={\boldsymbol 0}=(0,\ldots,0)$,
we fix
$\Delta^{(k)}({\boldsymbol 0})=1$.
As well as (\ref{eq:rTm}),
the symbol $A^i_j( k, l) $
denotes the
$ k \times  l$
 rectangular Toeplitz matrix 
 for the sequence
$\{h^i_j\}_{j=0}^\infty$
whose top left corner is $h^i_{j}$,
and
$h^i_j=0$ if $j<0$.
The second equality in (\ref{eq:bT2})
can be verified easily from
(\ref{eq:tenchi}). 
Henceforth
we suppose  
\[{\boldsymbol n}=(0,n, \ldots ,n)
\]
unless expressly stated otherwise.
We will often abbreviate
$\Delta^{(k)}({\boldsymbol n})$ for 
${\boldsymbol n}=(0,n, \ldots ,n)$
as $\Delta^{(k)}$.
Note that this convention is consistent with 
the description in
Remark~\ref{remark:bT}.
We also prepare the canonical basis 
$\{ {\boldsymbol e}_0,{\boldsymbol e}_1, \ldots, {\boldsymbol e}_{L-1} \}$ 
of ${\mathbb Z}^{L}$,
i.e. 
\[
{\boldsymbol e}_k=(0, \ldots,0,
\stackrel{ \stackrel{k}{\smile} }{1},0,\ldots,0).
\]

As seen in (\ref{eq:special}),
the coefficient 
\[
A(z)=\sum_{i=0}^{N+1} \frac{A_i}{z-u_i}
\quad 
(u_0=1, u_{N+1}=0)
\]
of
 the Fuchsian system ${\cal L}_{L,N}$ 
attached to the hypergeometric solution 
of ${\cal H}_{L,N}$
with $\kappa_0-\sum_{i=1}^N \theta_i=0$
(see Theorem~\ref{thm:hgsol})
is expressed as
$A_i= {\boldsymbol b}^{(i)} {\boldsymbol  c}^{(i)}$ $(0 \leq i\leq N)$,
where
\begin{align}
{\boldsymbol b}^{(0)} &=
h^0_0 \cdot
{}^{\rm T}
\left(0, \kappa_1,\kappa_2 ,
 \ldots, \kappa_{L-1} 
\right),
\quad
 {\boldsymbol c}^{(0)} =\left( 1, \frac{-1}{h^0_0}, \ldots, \frac{-1}{h^0_0} \right),
 \quad \text{and}
   \nonumber
 \\
 {\boldsymbol b}^{(i)} &=
- \theta_i 
 \cdot
 {}^{\rm T}
 \left(  1, {\ell_i}^{-1} (h_0^1),
 {\ell_i}^{-1} (h_0^2), 
 \ldots, 
  {\ell_i}^{-1} (h_0^{L-1})
 \right),
 \quad
  {\boldsymbol c}^{(i)} =( 1, 0, \ldots, 0 )
  \quad \text{for $1 \leq i \leq N$.}
   \label{eq:bchg}
\end{align}
Cf. (\ref{eq:resmat}).
Our next task is 
applying
the Schlesinger transformation to 
this Fuchsian system.

\subsection{Calculation of the Schlesinger transform (I)}
\label{subsect:cal1}

To derive the action of the Schlesinger transformation,
we need basically to deal with
(\ref{eq:AtoA}): 
\[
A
\mapsto \hat{A}
=
\sum_{i=0}^{N+1}
\frac{\hat{A}_i}{z-u_i}
=
R A R^{-1} + \frac{{\rm d} R}{{\rm d}z}R^{-1}.
\]
Namely, 
since
both $R=R(z)$ and $R^{-1}$ are polynomials in $z$,
each residue matrix 
$\hat{A}_i$
can be 
calculated by
\begin{equation}
\label{eq:hatai}
\hat{A}_i=
R(u_i) A_i R^{-1}(u_i)
\quad \text{for}
\quad 
0 \leq i \leq N+1.
\end{equation}
However,
thanks to (\ref{eq:Aid}),
it is rather easy to calculate the diagonal parts
even in the {\it general} case.
First we will demonstrate it.

The multiplier $R=R(z)$ of the Schlesinger transformation 
(see Theorem~\ref{thm:schl})
can be written, 
a little more specifically than (\ref{eq:R}),
as 
\[
R
= 
w^{-n}
\left(
\begin{bmatrix}
0 &  & &\\ 
*&Q_1^{(1)}(0)&  & \\ 
*&*&Q_2^{(2)}(0)& \\
\vdots&\vdots&\ddots&\ddots    &\\ 
*&*&\cdots&*&Q_{L-1}^{(L-1)}(0)
\end{bmatrix} 
+O(w)
\right);
\]
recall (\ref{eq:smatrix}).
Multiplying
\[
\Phi
=
\left(
\begin{bmatrix}
1 &  &\\  
 \vdots&\ddots& \\ 
*& \cdots&1
\end{bmatrix}
+O(w)
\right)
\Xi
\]
by $R$ from the left
yields 
$R \Phi = 
\hat{\Phi} 
\cdot {\rm diag \ } 
(w^{n(L-1)}, w^{-n}, \ldots,w^{-n})$
with
\[
\hat{\Phi}
=
\left(
\begin{bmatrix}
1 &  &\\  
 \vdots&\ddots& \\ 
*& \cdots&1
\end{bmatrix}
+O(w)
\right)
\hat{\Xi}
\]
as shown in Theorem~\ref{thm:schl},
where both
$\Xi$ and $\hat{\Xi}$ are diagonal matrices
independent of $w=1/z$.
We mention, 
without fear of repetition,
that
 the Hermite-Pad\'e approximation condition (\ref{eq:hpa})
assures 
\[
R  {\boldsymbol f}=O(w^{n (L-1)})
= 
\begin{bmatrix}
 \rho^0_{nL} 
 \\
  \rho^1_{nL} 
 \\
\vdots     \\
 \rho^{L-1}_{nL} 
 \end{bmatrix}
w^{n(L-1)} + 
(\text{terms of higher order})
\]
with
${\boldsymbol f}={}^{\rm T}(f_0, \ldots,f_{L-1})$
 denoting the first column of $\Phi$.
 We thus find the formula
\begin{align*}
 \hat{\Xi} {\Xi}^{-1}
&= {\rm diag \ } 
\left( \rho^0_{nL} , Q_1^{(1)}(0), Q_2^{(2)}(0),
\ldots,Q_{L-1}^{(L-1)}(0)
\right)
\\
&=
(-1)^n
\cdot
 {\rm diag \ } 
\left( (-1)^{nL} \frac{\Delta^{(L)}}{\Delta^{(1)}}, 
\frac{\Delta^{(1)}}{\Delta^{(2)}}, 
\frac{\Delta^{(2)}}{\Delta^{(3)}}, 
\ldots,
\frac{\Delta^{(L-1)}}{\Delta^{(L)}} 
\right)
\end{align*}
by virtue of Remark~\ref{remark:bT}.
Combining this with (\ref{eq:Aid})
under $x_i=1/u_i$:
\[
(A_i)_{\rm D}=x_i \frac{\partial}{\partial x_i}
\log \Xi,
\]
we arrive at the formulae
\begin{align} 
(\hat{A}_i)_{\rm D}-
(A_i)_{\rm D}
&=
x_i 
\frac{\partial}{\partial x_i}
\log 
{\rm diag \ } \left(
\frac{ \Delta^{(L)} }{ \Delta^{(1)} },
\frac{ \Delta^{(1)} }{ \Delta^{(2)} },
\ldots,
\frac{ \Delta^{(L-1)} }{ \Delta^{(L)} }
\right)
\nonumber
\\
&=
x_i
\
{\rm diag \ } \left(
 \frac{{\cal D}_i  \Delta^{(L)}  \cdot \Delta^{(1)} }{ \Delta^{(L)}  \Delta^{(1)} },
\frac{{\cal D}_i  \Delta^{(1)}  \cdot \Delta^{(2)} }{ \Delta^{(1)}  \Delta^{(2)} },
  \ldots,
 \frac{{\cal D}_i  \Delta^{(L-1)}  \cdot \Delta^{(L)} }{ \Delta^{(L-1)}  \Delta^{(L)} }
\right)
\label{eq:dhat}
\end{align}
for $1 \leq i \leq N$, where ${\cal D}_i$ denotes the Hirota differential with respect to $\partial/\partial x_i$.

Next we turn to the {\it particular} case, i.e.
the Fuchsian system ${\cal L}_{L,N}$ upon
the substitution (\ref{eq:special})
corresponding to the hypergeometric solution
(see Theorem~\ref{thm:hgsol}).
Write 
the residue matrices 
as 
$\hat{A}_i= \hat{\boldsymbol b}^{(i)} \hat{{\boldsymbol  c}}^{(i)}$ $(0 \leq i\leq N)$.
In order to reconstruct the canonical variables 
$(\hat{q}_k^{(i)},\hat{p}_k^{(i)})$ of the Schlesinger transform
of ${\cal H}_{L,N}$, 
it is only necessary to
know the quantities 
\begin{equation}
\label{eq:hatc}
\hat{{\boldsymbol  c}}^{(i)}
=
\left( \hat{c}_k^{(i)} \right)_{0 \leq k \le L-1}
=
\left(1, \hat{c}_1^{(i)}, \ldots ,\hat{c}_{L-1}^{(i)}\right)
 \quad \text{for} \quad  0 \leq i \leq N
\end{equation}
because we have already known 
from (\ref{eq:dhat})
 the diagonal entries
$ (\hat{A}_i)_{k,k}  =\hat{b}_k^{(i)} \hat{c}_k^{(i)}= -\hat{q}_k^{(i)} \hat{p}_k^{(i)}$
for $1 \leq i \leq N$;
cf. (\ref{eq:canvar}).
Calculations of (\ref{eq:hatc}) are a little complicated;
therefore, we will separately carry out 
the cases
$i=0$ and $i \neq 0$ in Sects.~\ref{subsect:cal2} and \ref{subsect:cal3},
respectively.
The result can be found in Sect.~\ref{subsect:result}.

\subsection{Calculation of the Schlesinger transform (II)}
\label{subsect:cal2}

Consider the row vector
${\boldsymbol c}^{(0)} R^{-1}(1)$
in view of (\ref{eq:hatai}) and $u_0=1$,
where
${\boldsymbol c}^{(0)}= (1,-1/h^0_0,\ldots,-1/h^0_0) $;
recall (\ref{eq:bchg}).

\paragraph{}
\underline{(i) The $0$th component}
\quad 
It follows from the determinantal representation of 
the polynomial
$P_0^{(0)}(w)$
and 
$P_l^{(0)}(w)
= 
\left[
 f_l P_0^{(0)}
\right]^{m-1}_{0}$
for $l > 0$
(see Proposition~\ref{prop:detp} and its sequel)
that
\[
P_l^{(0)}(1)
=
\frac{1}{ {\rm NP}^{(0)} }
\det 
\begin{bmatrix}
h_0^l +\cdots+h_{m-1}^l, 
h_0^l +\cdots+h_{m-2}^l, 
 \ldots,  
h_0^l +h_{1}^l, 
h_0^l,  
0 
\\
A^1_{m}(n,m+1) 
\\
\vdots
\\
A^{L-1}_{m}(n,m+1) 
\end{bmatrix}.
\]
Summation over $l=1,2, \ldots,L-1$ 
of this formula entails
\[
\sum_{l=1}^{L-1}
P_l^{(0)}(1)
=
\frac{1}{ {\rm NP}^{(0)} }
\det 
\begin{bmatrix}
h^0_0-h^0_m, 
h^0_0-h^0_{m-1}, 
 \ldots,  
h^0_0 -h^0_2, 
h^0_0-h^0_1,  
0 
\\
A^1_{m}(n,m+1) 
\\
\vdots
\\
A^{L-1}_{m}(n,m+1) 
\end{bmatrix}
\]
via the contiguity relation (\ref{eq:cont1}).
Hence the $0$th component of 
${\boldsymbol c}^{(0)} R^{-1}(1)$
reads
\begin{equation}
\label{eq:p01}
P^{(0)}_0(1)-\frac{1}{h^0_0} \sum_{l=1}^{L-1} P^{(0)}_l(1)
=
\frac{1}{h^0_0  {\rm NP}^{(0)} }
\det 
\begin{bmatrix}
h^0_m,h^0_{m-1}, \ldots, h^0_0
\\
A^1_{m}(n,m+1) 
\\
\vdots
\\
A^{L-1}_{m}(n,m+1) 
\end{bmatrix}
=
\frac{  \Delta^{(0)}({\boldsymbol n}+{\boldsymbol e}_0)}{h^0_0  \Delta^{(L)}({\boldsymbol n}) }.
\end{equation}
Here we have used (\ref{eq:np0}).

\paragraph{}
\noindent
\underline{(ii) The $k(> 0)$th component}
\quad
Similarly it holds for $l>0$
that
\[
P_l^{(k)}(1) 
=\frac{1}{{\rm NP}^{(k)}}
\det 
\begin{bmatrix}
h_0^l +\cdots+h_{m-1}^l, 
h_0^l +\cdots+h_{m-2}^l, 
 \ldots,  
h_0^l +h_{1}^l, 
h_0^l
\\
A^1_{m}(n,m) 
\\
\vdots
\\
A^{k-1}_{m}(n,m)
\\ 
A^{k}_{m}(n-1,m)
\\ 
A^{k+1}_{m-1}(n,m)
\\
\vdots \\
A^{L-1}_{m-1}(n,m) 
\end{bmatrix}.
\]
Taking a sum and using  (\ref{eq:cont1})
thus yield
\begin{align}
P^{(k)}_0(1)-\frac{1}{h^0_0} \sum_{l=1}^{L-1} P^{(k)}_l(1)
&=\frac{1}{ h^0_0 {\rm NP}^{(k)}}
\det 
\begin{bmatrix}
h^0_m,h^0_{m-1}, \ldots, h^0_1
\\
A^1_{m}(n,m) 
\\
\vdots
\\
A^{k-1}_{m}(n,m)
\\ 
A^{k}_{m}(n-1,m)
\\ 
A^{k+1}_{m-1}(n,m)
\\
\vdots \\
A^{L-1}_{m-1}(n,m) 
\end{bmatrix}
\nonumber
\\
&=\frac{(-1)^{nk+1}\Delta^{(k+1)}({\boldsymbol n}+{\boldsymbol e}_0-{\boldsymbol e}_k) }{ h^0_0\Delta^{(k)}  ({\boldsymbol n}) }.
\label{eq:pk1}
\end{align}
Here we have used (\ref{eq:npk}).

\paragraph{}
The ratio of
 (\ref{eq:p01}) and (\ref{eq:pk1})
 leads to the formula
\[
\hat{c}^{(0)}_k=
(-1)^{nk+1}
\frac{ \Delta^{(L)} ({\boldsymbol n})  
\Delta^{(k+1)} ({\boldsymbol n}+{\boldsymbol e}_0-{\boldsymbol e}_k) 
}{ \Delta^{(k)}  ({\boldsymbol n}) 
\Delta^{(0)}({\boldsymbol n}+{\boldsymbol e}_0) 
} .
\]

\subsection{Calculation of the Schlesinger transform (III)}
\label{subsect:cal3}

Consider
for $1 \leq i \leq N$
 the row vector
${\boldsymbol c}^{(i)} R^{-1}(u_i)$
in view of (\ref{eq:hatai}) and $u_i=1/x_i$,
which is nothing but 
the top row of the matrix 
$R^{-1}(u_i)$
due to
${\boldsymbol c}^{(i)}= (1,0,\ldots,0) $;
recall (\ref{eq:bchg}).

\paragraph{}
\underline{(i) The $0$th component}
\quad 
We are interested in the following determinant
\[
{x_i}^{-m}
P_0^{(0)}(x_i)
=
\frac{ {x_i}^{-m} }{ {\rm NP}^{(0)} }
\det 
\begin{bmatrix}
1,x_i, \ldots, {x_i}^m
\\
A^1_{m}(n,m+1) 
\\
\vdots
\\
A^{L-1}_{m}(n,m+1) 
\end{bmatrix}.
\]
Subtracting the $j$th column 
multiplied by $x_i$
from the $(j+1)$th column for 
$j=m,m-1,\ldots,1$ sequentially 
and using  the contiguity relation (\ref{eq:cont2}),
we thus obtain
\begin{align}
{x_i}^{-m}
P_0^{(0)}(x_i)
&=
\frac{ {x_i}^{-m} }{ {\rm NP}^{(0)} }
\det 
\begin{bmatrix}
1 &0, \ldots, 0
\\
A_m^1(n,1)&\ell_i \left( A^1_{m-1}(n,m) \right)
\\
\vdots&\vdots
\\
A_m^{L-1}(n,1)&\ell_i \left(A^{L-1}_{m-1}(n,m) \right)
\end{bmatrix}
=
\frac{ {x_i}^{-m} }{ {\rm NP}^{(0)} }
\ell_i 
\left(
\det 
\begin{bmatrix}
 A^1_{m-1}(n,m) 
\\
\vdots
\\
 A^{L-1}_{m-1}(n,m) 
\end{bmatrix}
\right)
\nonumber
\\
&=
\frac{(-1)^m {x_i}^{-m}  
\ell_i
\left(
\Delta^{(0)}({\boldsymbol n}) \right)}{ \Delta^{(L)}({\boldsymbol n}) }.
\label{eq:p0x}
\end{align}

\paragraph{}
\noindent
\underline{(ii) The $k(>0)$th component}
\quad 
In the same way as above, we find
\begin{align}
{x_i}^{-m+1}
P_0^{(k)}(x_i)
&=
\frac{ {x_i}^{-m+1} }{ {\rm NP}^{(k)} }
\det 
\begin{bmatrix}
1,x_i, \ldots, {x_i}^{m-1}
\\
A^1_{m}(n,m) 
\\
\vdots
\\
A^{k-1}_{m}(n,m)
\\ 
A^{k}_{m}(n-1,m)
\\ 
A^{k+1}_{m-1}(n,m)
\\
\vdots \\
A^{L-1}_{m-1}(n,m) 
\end{bmatrix}
=
\frac{ {x_i}^{-m+1} }{ {\rm NP}^{(k)} }
\ell_i 
\left(
\det 
\begin{bmatrix}
A^1_{m-1}(n,m-1) 
\\
\vdots
\\
A^{k-1}_{m-1}(n,m-1)
\\ 
A^{k}_{m-1}(n-1,m-1)
\\ 
A^{k+1}_{m-2}(n,m-1)
\\
\vdots \\
A^{L-1}_{m-2}(n,m-1) 
\end{bmatrix}
\right)
\nonumber
\\
&=
\frac{ (-1)^{n(L-k-1)} {x_i}^{-m+1} \ell_i \left(\Delta^{(k+1)}( {\boldsymbol n}-{\boldsymbol e}_k)  \right) }{ \Delta^{(k)} ( {\boldsymbol n}) }
\label{eq:pkx}
\end{align}

\paragraph{}
Hence we verify
from the ratio of
 (\ref{eq:p0x}) and (\ref{eq:pkx})
 that
\[
\hat{c}^{(i)}_k=
(-1)^{nk}
x_i
\frac{\Delta^{(L)} ({\boldsymbol n})
 \ell_i \left( \Delta^{(k+1)} ({\boldsymbol n}-{\boldsymbol e}_k) \right)
  }{\Delta^{(k)} ({\boldsymbol n}) 
  \ell_i \left( \Delta^{(0)} ({\boldsymbol n}) \right) } .
\]

\subsection{Result}
\label{subsect:result}

Summarizing the above we are led to the following result.

\begin{thm} 
\label{thm:hgi}
Let $n$ be a positive integer.
If $\kappa_0-\sum_{i=1}^N \theta_i=n(L-1)$
then
the Hamiltonian system
${\cal H}_{L,N}$ 
possesses 
a solution 
$(\hat{q}_k^{(i)},\hat{p}_k^{(i)})$ given by
\begin{align*}
\hat{q}_k^{(i)} &=
- x_i
\frac{ \Delta^{(0)}({\boldsymbol n}+{\boldsymbol e}_0)  
\ell_i 
\left( \Delta^{(k+1)} ({\boldsymbol n}-{\boldsymbol e}_k) \right)
}{ \Delta^{(k+1)} ({\boldsymbol n}+{\boldsymbol e}_0-{\boldsymbol e}_k) 
\ell_i \left( \Delta^{(0)} ({\boldsymbol n}) \right) 
 },
\\
\hat{q}_k^{(i)}  \hat{p}_k^{(i)}
&=
-x_i \frac{{\cal D}_i  \Delta^{(k)}({\boldsymbol n})  \cdot \Delta^{(k+1)}({\boldsymbol n})  }{ \Delta^{(k)}({\boldsymbol n})  \Delta^{(k+1)}({\boldsymbol n})  }
\end{align*}
under the correspondence
\[
\alpha_k=e_k-e_0, \quad 
\beta_i=  - \theta_i, \quad  
\gamma_k= e_k-e_0-\kappa_k+n
\]
of constant parameters,
where
${\boldsymbol n}=(0,n,\ldots,n) \in ({\mathbb Z}_{\geq 0})^L$.
\end{thm}

\begin{remark}\rm
We have in fact an alternative expression 
\[
\hat{q}_k^{(i)}  \hat{p}_k^{(i)}
=
\theta_i x_i 
\frac{
{\ell_i}^{-1} 
\left( \Delta^{(k)} ({\boldsymbol n}+{\boldsymbol e}_k) \right)
 \ell_i 
\left( \Delta^{(k+1)} ({\boldsymbol n}-{\boldsymbol e}_k) \right)
 }{ \Delta^{(k)}({\boldsymbol n})  \Delta^{(k+1)}({\boldsymbol n})  }
\]
of the above solution 
with no use of the Hirota differentials. 
This formula can be verified by calculating
$\hat{\boldsymbol b}^{(i)}$ 
in the same manner as Sects.~\ref{subsect:cal2} and \ref{subsect:cal3};
further details might be left to the reader.
\end{remark}

As seen in Theorem~\ref{thm:hgi},
we have constructed a particular solution of ${\cal H}_{L,N}$
expressed in terms of the block-Toeplitz determinant 
$\Delta^{(k)}({\boldsymbol n})$
whose entries are given by the hypergeometric functions.
Finally we shall rewrite $\Delta^{(k)}({\boldsymbol n})$
as an iterated hypergeometric integral:
if we remember that $h^i_j$ is defined as a moment
(see Sect.~\ref{subsect:prel}), 
then we observe that
\begin{align}
\nonumber
\Delta^{(k)}({\boldsymbol n})
&=
\int
\det \left[
\left[ {s_{0,j}}^{n_0+i-j} \right]_{
\begin{subarray}{l}
1 \leq i \leq |{\boldsymbol n}|\\
1 \leq j \leq n_0 
\end{subarray}},
\ldots,
\left[ {s_{k-1,j}}^{n_{k-1}+i-j} \right]_{
\begin{subarray}{l}
1 \leq i \leq |{\boldsymbol n}|\\
1 \leq j \leq n_{k-1} 
\end{subarray}},
\left[ {s_{k,j}}^{n_{k}+i-j-1} \right]_{
\begin{subarray}{l}
1 \leq i \leq |{\boldsymbol n}|\\
1 \leq j \leq n_{k} 
\end{subarray}},
\ldots,
\left[ {s_{L-1,j}}^{n_{L-1}+i-j-1} \right]_{
\begin{subarray}{l}
1 \leq i \leq |{\boldsymbol n}|\\
1 \leq j \leq n_{L-1} 
\end{subarray}}
\right]
\\
&
\nonumber
 \qquad
\times
\prod_{a=0}^{L-1}
\prod_{j=1}^{n_a} {\rm d}\mu_a(s_{a,j})
\\
\nonumber
&=
\int
\det \left[
\left[ {s_{0,j}}^{i-1} \right]_{
\begin{subarray}{l}
1 \leq i \leq |{\boldsymbol n}|\\
1 \leq j \leq n_0 
\end{subarray}},
\left[ {s_{1,j}}^{i-1} \right]_{
\begin{subarray}{l}
1 \leq i \leq |{\boldsymbol n}|\\
1 \leq j \leq n_1 
\end{subarray}},
\ldots,
\left[ {s_{L-1,j}}^{i-1}  \right]_{
\begin{subarray}{l}
1 \leq i \leq |{\boldsymbol n}|\\
1 \leq j \leq n_{L-1} 
\end{subarray}}
\right]
\prod_{a=0}^{k-1}
\prod_{j=1}^{n_a} 
{s_{a,j}}
\\
& 
\qquad 
\times
\prod_{a=0}^{L-1}
\prod_{j=1}^{n_a} 
{s_{a,j}}^{n_a-j}
{\rm d}\mu_a(s_{a,j})
\label{eq:intdelta}
\end{align}
through Fubini's theorem.
The Vandermonde determinant
shows that
\[
\det \left[
\left[ {s_{0,j}}^{i-1} \right]_{
\begin{subarray}{l}
1 \leq i \leq |{\boldsymbol n}|\\
1 \leq j \leq n_0 
\end{subarray}},
\left[ {s_{1,j}}^{i-1} \right]_{
\begin{subarray}{l}
1 \leq i \leq |{\boldsymbol n}|\\
1 \leq j \leq n_1 
\end{subarray}},
\ldots,
\left[ {s_{L-1,j}}^{i-1}  \right]_{
\begin{subarray}{l}
1 \leq i \leq |{\boldsymbol n}|\\
1 \leq j \leq n_{L-1} 
\end{subarray}}
\right]
=
\prod_{(a,b)>(c,d)} 
(s_{a,b}-s_{c,d}),
\]
where
$(a,b)>(c,d)$
is defined 
to mean
\[
\text{``$a>c$''}
\quad  \text{or} \quad
\text{``$a=c$ and $b >d$''.}
\]
Also it holds that
\[
\sum_{\sigma \in {\rm Sym \, }(n_a) } 
{\rm sgn \, } \sigma
\prod_{j=1}^{n_a} 
{s_{a,\sigma(j)}}^{n_a-j}
=\prod_{b<c} 
(s_{a,b}-s_{a,c}),
\]
where ${\rm Sym \, }(n_a)$
denotes the symmetric group on 
$\{1,2,\ldots,n_a\}$.
Therefore, 
because 
the value of integral (\ref{eq:intdelta})
 is 
invariant under a permutation of the variables $\{s_{a,b} \}_{1 \leq b \leq n_a}$
for each $0 \leq a \leq L-1$,
we conclude that
\[
\Delta^{(k)}({\boldsymbol n})
=
\prod_{a=0}^{L-1}
\frac{1}{n_a !}
\int
\left(\prod_{a=0}^{k-1}
\prod_{j=1}^{n_a} 
s_{a,j}
\right)
\left(
\prod_{(a,b)>(c,d)} 
(s_{a,b}-s_{c,d})
\right)
\prod_{a=0}^{L-1}
\left(
\prod_{b<c}
(s_{a,b}-s_{a,c})
\right)
\prod_{j=1}^{n_a} 
{\rm d}\mu_a(s_{a,j}).
\]

\appendix
\section{Verification of $(\ref{eq:bi})$}

Usually, we often normalize the Fuchsian system (\ref{eq:fs2}) 
so that $A_{N+2}=-\sum_{i=0}^{N+1} A_i $ 
(the residue matrix at $z=\infty$)
is diagonal
when considering its isomonodromic deformation;
see e.g. \cite{jmu81}.
But, in this paper, we adopt a different normalization
treating the two points 
$z=0$ and $z=\infty$ equally;
i.e.
we choose
$A_{N+1}$ (the residue matrix at $z=0$)
and $A_{N+2}$
to be upper and lower triangular,
respectively.
Note that the latter normalization is more versatile in
a general setting than the former.
In this appendix,
for a supplement to
Sect.~\ref{subsect:ss},
we demonstrate how to determine the coefficient $B_i=B_i(z)$ of 
the deformation equation (\ref{eq:defe}):
\[
\frac{\partial Y}{\partial u_i} =B_i Y
\]
 of (\ref{eq:fs2}).

Differentiating (\ref{eq:solatui}) 
with respect to $u_j$ tells us that
\begin{quote}
$\bullet$ if $i \neq j$ then
$\displaystyle\frac{\partial Y}{\partial u_j} Y^{-1}$ is holomorphic at
$z=u_i$;
\\
$\bullet$ if $i=j$ then 
$\displaystyle
 \frac{\partial Y}{\partial u_i} Y^{-1}=
- \frac{A_i}{z-u_i}
+(\text{holomorphic at $z=u_i$}).$
\end{quote}
Similarly, we observe that
\begin{equation}
\label{eq:solatzero}
 \frac{\partial Y}{\partial u_i} Y^{-1}=\frac{\partial \Psi}{\partial u_i} \Psi^{-1}
 = 
 \begin{bmatrix}
 0 & * & \cdots & * \\
    & 0 & \ddots & \vdots \\
    &    &  \ddots & * \\
    &    &             & 0
 \end{bmatrix} +O(z)
 \quad \text{near $z=0$}
\end{equation}
and
\begin{equation}
\label{eq:solatinfty}
 \frac{\partial Y}{\partial u_i} Y^{-1}=\frac{\partial \Phi}{\partial u_i} \Phi^{-1}
 = 
 \begin{bmatrix}
 * &  & \\
 \vdots   & \ddots &  \\
*    &  \cdots  &   * 
 \end{bmatrix} +O(w)
 \quad \text{near $z=1/w=\infty$}.
\end{equation}
Here we have used the assumption that
the connection matrix $C$ does not depend on $u_i$.
Consequently, 
$\displaystyle B_i= \frac{\partial Y}{\partial u_i} Y^{-1}$
is a rational function matrix in $z$
that has only a simple pole at $z=u_i$ with
residue $-A_i$
and
thus
\[
B_i= \frac{A_i}{u_i-z}+K_i,
\]
where
$K_i$ is a constant matrix.
Substituting $z=\infty$ ($w=0$)
in (\ref{eq:solatinfty}) shows that $K_i$ is a lower triangular matrix. 
Therefore, substituting $z=0$ in (\ref{eq:solatzero}),
we  conclude that
$K_i=- (A_i)_{\rm LT}/u_i$.

\small
\paragraph{\it Acknowledgement.}
The authors are deeply grateful to
Shuhei Kamioka  
for giving them an exposition on various multi-dimensional continued fractions.
They appreciate Satoshi Tsujimoto 
his kind information about literature on rational approximations.
Also, they have benefited from invaluable discussions with Yasuhiko Yamada.
This work was supported
in part by a grant-in-aid from the Japan Society for the Promotion of Science (Grant Number: 25800082 and 25870234).


\end{document}